\documentclass[a4paper]{amsart}

\usepackage{amssymb}
\usepackage{graphics}
\usepackage{color}
\usepackage[all]{xy}
\CompileMatrices

\definecolor{changes}{rgb}{0.1,0.65,0.03}

\newtheorem*{introthm}{Theorem}
\newtheorem{theorem}{Theorem}[section]
\newtheorem{proposition}[theorem]{Proposition}
\newtheorem{lemma}[theorem]{Lemma}
\newtheorem{corollary}[theorem]{Corollary}

\theoremstyle{definition}
\newtheorem*{introex}{Example}
\newtheorem{definition}[theorem]{Definition}
\newtheorem{remark}[theorem]{Remark}
\newtheorem{example}[theorem]{Example}

\theoremstyle{remark}
\newtheorem{step}{Step}
\newtheorem*{ecase}{Elliptic case}
\newtheorem*{pcase}{Parabolic case}
\newtheorem*{hcase}{Hyperbolic case}

\numberwithin{equation}{section}

\def\Rho{{\rm R}}
\def\gittera{N} 
\def\dualgittera{M} 
\def\unterdualgittera{L} 
\def\linabba{F} 

\def\kegela{\sigma} 
\def\dualkegela{\omega} 
\def\polytopa{\Pi} 
\def\polyedera{\Delta} 
\def\polyederb{\Gamma} 
\def\relint#1{{\rm relint}(#1)} 

\def\gewkegel{\omega} 
\def\gkfaecher{\Lambda} 
\def\gkfkegela{\lambda} 
\def\gkfkegelb{\gamma} 

\def\eval{{\rm eval}} 
\DeclareMathOperator{\lin}{lin} 

\def\torusa{T} 
\def\quotabb{\pi} 
\def\morphX{\varphi} 
\def\morphY{\psi} 

\DeclareMathOperator{\PPDiv}{PPDiv} 
\DeclareMathOperator{\CPL}{CPL} 
\DeclareMathOperator{\ddiv}{div} 

\def\poldiva{\mathfrak{D}} 
\def\cpla{\mathfrak{h}} 
\def\polfunca{\mathfrak{f}} 

\def\quot{/\!\!/} 
\def\mal{\! \cdot \!} 
\def\til#1{\widetilde{#1}} 
\def\blk#1{\overline{#1}} 
\def\blkk#1{\bar{#1}} 
\def\kst{; \; } 
\def\rd#1{\lfloor \! #1 \! \rfloor} 

\newcommand{\T}{{T}}
\newcommand{\Ta}{\T}

\DeclareMathOperator{\CDiv}{CaDiv}
\DeclareMathOperator{\WDiv}{Div}
\DeclareMathOperator{\Pic}{Pic}
\DeclareMathOperator{\Pol}{Pol}

\newcommand{\funcX}{\mathfrak X}
\newcommand{\funcD}{\blk{\poldiva}}

\newcommand{\N}{\mathbb N}

\newcommand{\Z}{\mathbb Z}
\newcommand{\PP}{\mathbb P}
\newcommand{\C}{\mathbb C}
\newcommand{\K}{\mathbb K}
\newcommand{\Q}{\mathbb Q}

\DeclareMathOperator{\Spec}{Spec}
\DeclareMathOperator{\Proj}{Proj}
\DeclareMathOperator{\Hom}{Hom}

\DeclareMathOperator{\cone}{cone}
\DeclareMathOperator{\Div}{Div}
\DeclareMathOperator{\kdiv}{div}
\DeclareMathOperator{\Supp}{Supp} 

\newcommand{\kprojlim}{\mathop{\varprojlim}\limits} 
\newcommand{\kdirlim}{\mathop{\varinjlim}\limits} 
\DeclareMathOperator{\id}{id}

\begin{document}

\title[Polyhedral Divisors]
      {Polyhedral Divisors and \\ Algebraic Torus Actions}  

\author[K.~Altmann]{Klaus Altmann$^{1}$}
\thanks{$^{1}$partially supported by MSRI Berkeley, 
        CA and SFB 647 of the DFG}
\address{Fachbereich Mathematik und Informatik, 
Freie Universit\"at Berlin,
Arnimalle 3, 
14195 Berlin, 
Germany}
\email{altmann@math.fu-berlin.de}
\author[J.~Hausen]{J\"urgen Hausen$^{2}$} 
\thanks{$^{2}$partially supported by SWP 1094 of the DFG}
\address{Mathematisches Institut, 
         Universit\"at T\"ubingen, 
         Auf der Morgenstelle 10,
         72076 T\"ubingen, Germany}
\email{hausen@mail.mathematik.uni-tuebingen.de}


\begin{abstract}
We provide a complete description of normal affine varieties
with effective algebraic torus action in terms of what we call
proper polyhedral divisors on semiprojective varieties.
Our approach extends classical cone constructions of Dolgachev, 
Demazure and Pinkham to the multigraded case,
and it comprises the theory of affine toric varieties.
\end{abstract}

\maketitle

\section*{Introduction}
We present a complete description of $n$-dimensional,
normal, affine varieties with an effective action of
a $k$-dimensional algebraic torus in terms of
``proper polyhedral divisors'' living on semiprojective varieties
of dimension $n-k$.
Our approach comprises two well known theories:
on the one hand, for varieties with an almost transitive 
torus action ($k=n$), our description specializes to the theory
of affine toric varieties~\cite{Fulton}, 
and on the other, for $\C^{*}$-actions ($k=1$), we recover
classical constructions of generalized affine cones of 
Dolgachev~\cite{Dol}, Demazure~\cite{Dem} and 
Pinkham~\cite{Pi}. 

Besides the special cases $k=1$ and $k=n$, 
also the case $k=n-1$ is studied
by other authors, even for not necessarily affine
$T$-varieties $X$.
In the last chapter of~\cite{Ke}, the combinatorial
methods for toroidal varieties developed in this book
are applied to study torus actions of codimension one;
see~\cite{Vollmert} for a comparison of this approach
and ours.
In~\cite{Timashev}, Timashev presents a general theory
of reductive group actions of complexity one on normal
algebraic varieties.
Specializing his language of hypercones to the case 
of a torus action,
he obtains a picture quite similar to ours, see
Example~4.1 in loc.~cit.
Moreover, there is recent work by Flenner and Zaidenberg
on affine $\K^*$-surfaces~\cite{FlZa}, which fits into our
framework, see Example~\ref{XXX}.
Finally, the analogous setting is also studied in symplectic
geometry, see for example the treatment in \cite{KarshonTolman}
using the moment map.

Let us outline the main results of the 
present paper. 
Let $Y$ be a normal semiprojective
variety, where ``semiprojective'' merely 
means that $Y$ is 
projective over some affine variety.
In order to introduce the notion of 
a proper polyhedral divisor on $Y$, 
consider a linear combination
\begin{eqnarray*}
\poldiva
& = & 
\sum 
\polyedera_i \otimes D_i 
\end{eqnarray*}
where the $D_i$ are prime divisors on 
$Y$, the coefficients $\Delta_i$ 
are convex polyhedra in a rational vector
space $N_\Q = \Q \otimes N$ with a free 
finitely generated abelian group $N$, 
and all $\Delta_i$ have a common pointed
cone $\sigma \subset N_\Q$ as their
tail cone
(see Section~\ref{section1} for the 
precise definitions).

Let $M := \Hom(N,\Z)$ be the dual 
of $N$, and write $\sigma^\vee \subset M_\Q$ 
for the dual cone. 
Then the above $\poldiva$
defines an evaluation map
into the group of rational 
Weil divisors on $Y$: 
$$
\sigma^\vee 
\; \to \; 
\WDiv(Y),
\qquad 
u  \; \mapsto \; \mathfrak{D}(u) 
:= 
\sum \min_{v \in \Delta_i} \langle u,v \rangle D_i. 
$$
We say that $\poldiva$ is a 
{\em proper polyhedral divisor\/}
if any evaluation $\poldiva(u)$ 
is a semiample rational Cartier divisor,
being big whenever $u$ belongs to the 
relative interior of the cone 
$\sigma^\vee$.

The evaluation map $u \mapsto \poldiva(u)$
turns out to be piecewise linear and 
convex in the sense that 
the difference $\poldiva(u+u') - \big(\poldiva(u) +\poldiva(u')\big)$ 
is always effective.
This convexity property enables us
to define a graded algebra of global sections:
\begin{eqnarray*}
A 
& := &
\bigoplus_{u \in \sigma^\vee \cap M} \Gamma(Y,\mathcal{O}(\poldiva(u))).
\end{eqnarray*}
As we will prove, among other things, in Theorem~\ref{data2variety},
this ring is normal and finitely generated. 
Thus, it gives rise to a normal affine variety $X := \Spec(A)$, 
and the $M$-grading of $A$ defines an effective action of the 
torus $T := \Spec(\C[M])$ on~$X$.

\begin{introex}
Let $Y = \PP^1$ and $N = \Z^2$. 
The vectors $(1,0)$ and $(1,12)$
generate a pointed convex cone 
$\sigma$ in $N_\Q = \Q^2$,
and we consider the polyhedra
$$
\Delta_0 \; = \; 
\left(\frac{1}{3},0\right)+\sigma,
\qquad
\Delta_1 \; = \; \left(-\frac{1}{4},0\right) + \sigma,
\qquad
\Delta_\infty \; = \; (\{0\}\times[0,1])+\sigma.
$$
Attaching these polyhedra as coefficients to the 
points $0, 1, \infty$ on the projective line, we 
obtain a proper polyhedral divisor
\begin{eqnarray*}
\poldiva^{E_6}
& = &
\Delta_0 \otimes \{0\} 
\; + \; 
\Delta_1 \otimes \{1\}
\; + \;
\Delta_\infty \otimes \{\infty\}. 
\end{eqnarray*}
In this situation, we may even represent our proper polyhedral 
divisor by a little picture as follows: 
\medskip

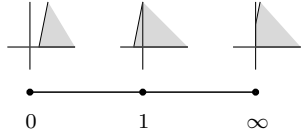
\vbox{
\begin{center}
\begin{picture}(0,0)%
\includegraphics{e6poldiv.pstex}%
\end{picture}%
\setlength{\unitlength}{1243sp}%
\begingroup\makeatletter\ifx\SetFigFont\undefined%
\gdef\SetFigFont#1#2#3#4#5{%
  \reset@font\fontsize{#1}{#2pt}%
  \fontfamily{#3}\fontseries{#4}\fontshape{#5}%
  \selectfont}%
\fi\endgroup%
\begin{picture}(5884,2650)(439,-2239)
\put(3061,-2131){\makebox(0,0)[lb]{\smash{\SetFigFont{8}{9.6}{\familydefault}{\mddefault}{\updefault}{\color[rgb]{0,0,0}$1$}%
}}}
\put(811,-2131){\makebox(0,0)[lb]{\smash{\SetFigFont{8}{9.6}{\familydefault}{\mddefault}{\updefault}{\color[rgb]{0,0,0}$0$}%
}}}
\put(5221,-2131){\makebox(0,0)[lb]{\smash{\SetFigFont{8}{9.6}{\familydefault}{\mddefault}{\updefault}{\color[rgb]{0,0,0}$\infty$}%
}}}
\end{picture}

\medskip

\end{center}
}
\noindent
As we shall see in Section~\ref{sec:examples}, 
the proper polyhedral divisor
$\poldiva^{E_6}$ describes the affine threefold
$X = V(z_1^3+z_2^4+z_3z_4)$ in $\C^4$ with the action
of $T = (\C^*)^2$ given by 
\begin{eqnarray*}
t \mal z 
& = & 
(t_1^4z_1, t_1^3z_2, t_2z_3, t_1^{12}t_2^{-1}z_4).
\end{eqnarray*}
\end{introex}

Assigning to the pp-divisor $\poldiva$ the affine 
$T$-variety $X$, as indicated, turns out to be functorial.
Moreover, a canonical construction, based on the chamber
structure of the set of GIT-quotients of $X$, shows
that in fact every normal affine variety with effective 
torus action arises from a proper polyhedral divisor.
These results can be summarized as follows,
see Theorems~\ref{variety2data} and
Proposition~\ref{prop:functorial}.

\begin{introthm}
The assignment $\poldiva \mapsto X$ defines
an essentially surjective faithful covariant functor 
from the category of proper polyhedral divisors on 
semiprojective varieties to the category of normal 
affine varieties with effective torus action.
\end{introthm}

After localizing the category of proper polyhedral divisors 
by the maps coming from
(birational) modifications of the semiprojective base
varieties, 
we even arrive at an equivalence of categories,
compare Corollary~\ref{equivcat}.
In particular, these results allow the determination of when two
proper polyhedral divisors define (equivariantly) isomorphic
varieties.

As an application, in Section~\ref{sec:orbits},
we provide 
a description of the collection of $T$-orbits 
of an affine $T$-variety
$X$ in terms of its defining pp-divisor $\mathfrak{D}$, 
and we indicate how to read local orbit data from 
$\mathfrak{D}$. 
Moreover, we indicate in Section~\ref{sec:examples}
a recipe for the computation of the pp-divisor
of a given affine variety with torus action.

In a subsequent paper, we will deal with {\em non-affine}
$T$-varieties $X$. Then, the coefficients of the former
polyhedral divisors on $Y$ will turn into polyhedral complexes.

\tableofcontents

We would like to especially thank 
J.~A.~Christophersen for valuable and stimulating discussions.
Moreover, we are grateful to the referee,
as well as to I.V.~Arzhantsev and R.~Vollmert 
for helpful remarks on earlier versions 
of the manuscript.

\section{Tailed polyhedra}

In this section, we introduce and discuss
the groups of tailed polyhedra,
which will serve later as the group of 
coefficients for our polyhedral divisors.
While setting 
the definitions and 
statements, we also fix our notation 
from convex geometry, and we recall some 
basic facts needed later.
For further background on convex geometry 
we refer to standard text books, 
like~\cite{Gr} and~\cite{Ro}. 

{From} here on, $\gittera$ denotes a
{\em lattice}, i.e.\ a finitely generated free 
abelian group.
The rational vector space associated to $\gittera$
is denoted by $\gittera_{\Q} := \Q \otimes_{\Z}  \gittera$,
and, given a linear map 
$\linabba \colon \gittera \to \gittera'$, 
we also write 
$\linabba \colon \gittera_{\Q} \to \gittera'_{\Q}$
for the induced map 
of rational vector spaces.
The dual lattice of $\gittera$ is
$\dualgittera := \Hom(\gittera,\Z)$,
and we denote the corresponding 
pairing by
$$ 
\dualgittera \times \gittera \; \to \; \Z,
\qquad
(u,v) \; \mapsto \; \langle u, v \rangle. 
$$

By a {\em polyhedron\/} in $\gittera_{\Q}$,
we mean a convex polyhedron, i.e.\ 
the intersection of finitely many closed 
affine half spaces in $\gittera_{\Q}$.
If $\polyedera'$ is a face of a 
polyhedron $\polyedera$ in $N_{\Q}$, 
then we write
$\polyedera' \preceq \polyedera$.
For a polyhedron $\polyedera$ in 
$\gittera_{\Q}$, 
we denote by $\relint{\polyedera}$ 
its {\em relative interior},
i.e.\ the set obtained by removing all
proper faces from $\polyedera$.

For us, a {\em cone\/} in $\gittera_{\Q}$
is always a convex, polyhedral cone,
i.e.\ the intersection of finitely 
many closed linear half spaces in 
$\gittera_{\Q}$.
The {\em dual cone\/} $\sigma^{\vee}$ 
of a cone $\kegela$ in
$\gittera_\Q$ lives in the dual 
vector space $\dualgittera_\Q$
and consists of all linear forms of 
$\dualgittera_{\Q}$ that are 
nonnegative along $\kegela$.
A cone is {\em pointed\/} if it does not
contain any line.

The set of all polyhedra in $\gittera_{\Q}$
comes with a natural abelian semigroup
structure: one defines the 
{\em Minkowski sum\/} of two polyhedra
$\polyedera_1$ and $\polyedera_2$ in
$\gittera_{\Q}$ to be the polyhedron 
\begin{eqnarray*}
\polyedera_1 + \polyedera_2
& = & 
\{v_1+v_2; \; v_i \in \polyedera_i\}.
\end{eqnarray*}

Any polyhedron $\polyedera$ in $\gittera_{\Q}$ 
allows a Minkowski sum decomposition 
$\polyedera = \polytopa + \kegela$
where $\polytopa \subset \gittera_{\Q}$
is a {\em polytope}, i.e.\ the convex hull of
finitely many points, and 
$\sigma \subset \gittera_{\Q}$ 
is a cone. 
In this decomposition, the 
{\em tail cone\/} $\sigma$ is unique; 
in the literature it is also called the 
{\em recession cone\/} of
$\polyedera$ and is given by
\begin{eqnarray*}
\sigma
& = &
\{v \in \gittera_{\Q}; \; 
  v' + tv \in \polyedera
  \text{ for all }
  v' \in \Delta, \, t \in \Q_{\ge 0}
\}.
\end{eqnarray*}

\begin{definition}
Let $\kegela$ be a pointed cone in $\gittera_{\Q}$.
\begin{enumerate}
\item 
By a {\em $\kegela$-tailed polyhedron\/} 
(or {\em $\kegela$-polyhedron}, in short) 
in $\gittera_{\Q}$,
we mean a polyhedron $\polyedera$ 
in $\gittera_{\Q}$ 
having the cone $\sigma$ as its tail cone.
We denote the set of all $\kegela$-polyhedra
in $\gittera_{\Q}$ by 
$\Pol^{+}_{\kegela}(\gittera_{\Q})$.
\item
We call $\polyedera \in \Pol^{+}_{\kegela}(\gittera_{\Q})$
{\em integral\/} if 
$\polyedera = \polytopa + \kegela$
holds with a polytope 
$\polytopa  \subset \gittera_{\Q}$
having its vertices in~$\gittera$.
We denote the set of all integral 
$\kegela$-polyhedra
in $\gittera_{\Q}$ by 
$\Pol^{+}_{\kegela}(\gittera)$.
\end{enumerate}
\end{definition}

Note that the Minkowski sum 
$\polyedera_1 + \polyedera_2$ 
of two $\kegela$-polyhedra 
$\polyedera_1$ and $\polyedera_2$ in
$\gittera_{\Q}$ is 
again a $\kegela$-polyhedron
in $\gittera_{\Q}$.
Thus, together with Minkowski addition, 
$\Pol^{+}_{\kegela}(\gittera_{\Q})$
is an abelian monoid; its neutral element is
$\kegela \in \Pol^{+}_{\kegela}(\gittera_{\Q})$, 
and $\Pol^{+}_{\kegela}(\gittera)
\subset \Pol^{+}_{\kegela}(\gittera_{\Q})$
is a submonoid.

\medskip

\vbox{
\begin{center}
\begin{picture}(0,0)%
\includegraphics{minkowski.pstex}%
\end{picture}%
\setlength{\unitlength}{1243sp}%
\begingroup\makeatletter\ifx\SetFigFont\undefined%
\gdef\SetFigFont#1#2#3#4#5{%
  \reset@font\fontsize{#1}{#2pt}%
  \fontfamily{#3}\fontseries{#4}\fontshape{#5}%
  \selectfont}%
\fi\endgroup%
\begin{picture}(11352,3706)(361,-3305)
\put(4501,-1636){\makebox(0,0)[lb]{\smash{\SetFigFont{8}{9.6}{\familydefault}{\mddefault}{\updefault}{\color[rgb]{0,0,0}$\Pi_2$}%
}}}
\put(6301,-961){\makebox(0,0)[lb]{\smash{\SetFigFont{8}{9.6}{\familydefault}{\mddefault}{\updefault}{\color[rgb]{0,0,0}$\sigma$}%
}}}
\put(10801,-961){\makebox(0,0)[lb]{\smash{\SetFigFont{8}{9.6}{\familydefault}{\mddefault}{\updefault}{\color[rgb]{0,0,0}$\sigma$}%
}}}
\put(1576,-3211){\makebox(0,0)[lb]{\smash{\SetFigFont{8}{9.6}{\familydefault}{\mddefault}{\updefault}{\color[rgb]{0,0,0}$\Delta_1$}%
}}}
\put(9676,-3211){\makebox(0,0)[lb]{\smash{\SetFigFont{8}{9.6}{\familydefault}{\mddefault}{\updefault}{\color[rgb]{0,0,0}$\Delta_1+\Delta_2$}%
}}}
\put(2251,-961){\makebox(0,0)[lb]{\smash{\SetFigFont{8}{9.6}{\familydefault}{\mddefault}{\updefault}{\color[rgb]{0,0,0}$\sigma$}%
}}}
\put(5851,-3211){\makebox(0,0)[lb]{\smash{\SetFigFont{8}{9.6}{\familydefault}{\mddefault}{\updefault}{\color[rgb]{0,0,0}$\Delta_2$}%
}}}
\put(361,-511){\makebox(0,0)[lb]{\smash{\SetFigFont{8}{9.6}{\familydefault}{\mddefault}{\updefault}{\color[rgb]{0,0,0}$\Pi_1$}%
}}}
\end{picture}

\bigskip

{\small
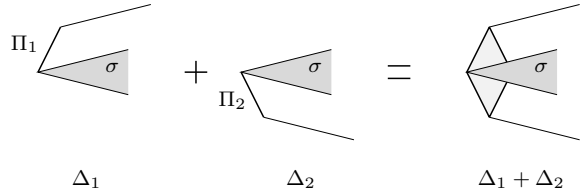
Minkowski addition of two 
$\sigma$-polyhedra 
$\polyedera_i = \Pi_i + \sigma$ in~$\Q^2$.
}
\end{center}
}

%

\goodbreak

\begin{definition}
Let $\kegela$ be a pointed cone in $\gittera_{\Q}$.
\begin{enumerate}
\item
The {\em group of $\sigma$-polyhedra\/}
is the Grothendieck group of 
$\Pol^{+}_{\kegela}(\gittera_{\Q})$;
we denote it by $\Pol_{\kegela}(\gittera_{\Q})$.
\item
The {\em group of integral $\sigma$-polyhedra\/}
is the Grothendieck group of 
$\Pol^{+}_{\kegela}(\gittera)$;
we denote it by $\Pol_{\kegela}(\gittera)$.
\end{enumerate}
\end{definition}

The key to basic properties
of these groups is a version 
of the general correspondence 
between convex sets 
and so-called support 
functions~\cite[Theorem~13.2]{Ro}
adapted to the setting of tailed 
polyhedra.  
In order to state this adapated version,
we firstly have to recall further 
notions from convex geometry.

A {\em quasifan\/} $\gkfaecher$ in $\dualgittera_{\Q}$ 
is a finite collection of cones in $\dualgittera_{\Q}$ 
with the following properties, compare~\cite[1.2]{Reid0}:
for any $\gkfkegela\in \gkfaecher$,
all the faces  
$\gkfkegela'\preceq\gkfkegela$ belong to $\gkfaecher$,
and, for any two $\gkfkegela_{i}\in\gkfaecher$, the
intersection $\gkfkegela_{1}\cap\gkfkegela_{2}$ is
a face of each $\gkfkegela_{i}$.
The {\em support\/} of a quasifan is the union of its cones.
A quasifan is called a {\em fan\/} if all its cones
are pointed.

To every polyhedron $\polyedera$ in $\gittera_{\Q}$,
one associates its {\em normal quasifan\/}
$\gkfaecher(\polyedera)$ in $\dualgittera_{\Q}$;
the faces $F \preceq \polyedera$
are in order reversing bijection with 
the cones of $\gkfaecher(\polyedera)$
via
\begin{eqnarray*}
F 
& \mapsto &
\lambda(F)
\; := \; 
\{u \in \dualgittera_{\Q}; \; 
\langle u, v - v' \rangle \ge 0
\text{ for all } 
v \in \polyedera, \, v' \in F\}.
\end{eqnarray*}
%
It is a basic observation that the
normal quasifan 
$\gkfaecher(\polyedera_1 + \polyedera_2)$
of a Minkowski sum is 
supported on the intersection of the supports of
the normal quasifans $\gkfaecher(\polyedera_1)$
and $\gkfaecher(\polyedera_2)$
and, moreover, equals the coarsest common refinement of both.

\begin{lemma}
\label{normalfansupport}
Let $\kegela \subset \gittera_{\Q}$ 
be a pointed cone,
and let 
$\polyedera \in \Pol^{+}_{\kegela}(\gittera_{\Q})$.
Then the normal quasifan 
$\gkfaecher(\polyedera)$ has
the dual cone 
$\kegela^{\vee} \subset \dualgittera_{\Q}$ 
as its support.
\end{lemma}

\begin{proof}
For every face $F \preceq \polyedera$,
the set
$\{v - v'; \; 
   v \in \polyedera, \, 
   v' \in F\}$
contains the tail cone $\kegela$.
Dualizing yields that the cone of
$\gkfaecher(\polyedera)$ corresponding
to $F$ is
contained in $\kegela^\vee$.
Conversely, every $u \in \kegela^\vee$
attains its minimum along some face of
$F \preceq \polyedera$,
and hence belongs to
a cone of $\gkfaecher(\polyedera)$.
\end{proof}

Since we require $\kegela$ to be pointed, 
$\kegela^\vee$ is of full dimension. 
The Lemma thus implies that for any 
$\sigma$-polyhedron $\polyedera$, 
the maximal cones of 
$\gkfaecher(\polyedera)$ are of full 
dimension, and hence, the minimal 
faces of $\polyedera$ are vertices, 
i.e.\ are of dimension zero. 
The vertices of $\polyedera$ are 
vertices of any polytope $\Pi$
with $\polyedera = \Pi + \kegela$,
and we may canonically write 
$\polyedera = \Pi_0 + \sigma$,
where $\Pi_0$ is the convex hull
of the vertices of $\polyedera$.

Next, we have to recall the definition
of the {\em support function\/} 
associated with a convex set 
$\polyedera$ in $\gittera_{\Q}$;
this is the map given by
$$
h_{\polyedera}
\colon 
M_{\Q}
\; \to \; 
\Q \cup\{-\infty\},
\qquad
u \; \mapsto \; 
\inf_{v \in \Delta} \langle u, v \rangle.
$$
The {\em domain\/} of this function is 
the subset of $\dualgittera_{\Q}$
where it takes values in $\Q$. 
Here are the basic properties of
the support function 
of a $\kegela$-polyhedron.

\begin{lemma}
\label{poly2fct}
Let $\kegela$ be a pointed cone in $\gittera_{\Q}$, 
let $\polyedera \in \Pol^{+}_{\kegela}(\gittera_{\Q})$,
and let $h_{\polyedera}$ be the corresponding
support function.
\begin{enumerate}
\item 
The function $h_{\Delta}$ has the dual 
cone $\kegela^{\vee}$ as its domain,
and it is linear on each cone of the 
normal quasifan $\gkfaecher(\polyedera)$.
\item
The function $h_{\Delta}$ is 
convex, that means that for any two 
vectors $u_1,u_2 \in \kegela^{\vee}$ we have 
\begin{eqnarray*}
h_{\polyedera}(u_1)
+
h_{\polyedera}(u_2)
& \le & 
h_{\polyedera}(u_1+u_2).
\end{eqnarray*}
Moreover, strict inequality holds if and only
if the vectors $u_1,u_2 \in \kegela^{\vee}$ do not
belong to the same maximal cone of $\Lambda(\Delta)$.
\end{enumerate}
\end{lemma}

\begin{proof}
The statements are standard in the case that
$\polyedera$ is a polytope of full dimension,
see~\cite[Appendix~A]{Od}; the simple proofs 
given there are easily adapted to our setting. 
\end{proof}

As usual, we say that a function 
$h \colon M_{\Q} \to \Q \cup \{-\infty\}$  
with a cone $\omega \subset M_{\Q}$
as its domain is {\em piecewise linear\/}
if there is a quasifan $\gkfaecher$ 
having $\omega$
as its support such that $h$ is linear on the cones 
of $\gkfaecher$.
We denote the set of convex piecewise linear functions 
on $M_{\Q}$ 
having a given cone $\omega$ as its domain
by $\CPL_{\Q}(\omega)$. 
Together with pointwise 
addition, $\CPL_{\Q}(\omega)$
is an abelian monoid.

\begin{proposition}
\label{poly2suppfunct}
Let $\kegela$ be a pointed cone in $N_{\Q}$.
Then the map 
$
\Pol^{+}_{\kegela}(\gittera_{\Q}) 
\to 
\CPL_{\Q}(\kegela^{\vee})$,
$\polyedera \mapsto h_{\polyedera}$ is 
an isomorphism of abelian semigroups.
\end{proposition}

\begin{proof}
According to Lemma~\ref{poly2fct}, 
the map is well-defined, and it is
easily checked to be a monoid 
homomorphism.
Moreover, the assignment 
\begin{eqnarray*}
h 
& \mapsto & 
\polyedera_h :=
\{
v \in \gittera_\Q; \;
\langle u, v \rangle \ge h(u)
\text{ for all }
u \in \kegela^{\vee}
\}
\end{eqnarray*}
associates to any $h \in \CPL_{\Q}(\kegela^{\vee})$
a $\kegela$-polyhedron, and it is directly checked
that this gives the inverse homomorphism.
\end{proof}

As announced, we now apply this observation to 
provide basic properties of the groups of 
$\sigma$-polyhedra. 
 
\begin{proposition}
\label{cancel}
Let $\kegela$ be a pointed cone in $\gittera_\Q$.
Then $\Pol^{+}_{\kegela}(\gittera_{\Q})$
and $\Pol^{+}_{\kegela}(\gittera)$ 
are abelian monoids with cancellation
law. Their respective groups of units are 
$$
\Pol^{+}_{\kegela}(\gittera_{\Q})^* 
\; = \:
\{v + \sigma; \; v \in \gittera_{\Q}\},
\qquad
\Pol^{+}_{\kegela}(\gittera)^* 
\; = \:
\{v + \sigma; \; v \in \gittera\}.
$$
\end{proposition}

\begin{proof}
Clearly, $\CPL_{\Q}(\kegela^{\vee})$ is 
an abelian monoid with cancellation law.
By Proposition~\ref{poly2suppfunct}, 
the same holds for 
$\Pol^{+}_{\kegela}(\gittera_{\Q})$
and the submonoid 
$\Pol^{+}_{\kegela}(\gittera) 
\subset
\Pol^{+}_{\kegela}(\gittera_{\Q})$. 
Moreover, the polyhedra $v + \kegela$ 
correspond to the linear functions 
$u \mapsto \langle u,v \rangle$, which 
are invertible in $\CPL_{\Q}(\kegela^{\vee})$.
Since the negative of a nonlinear 
convex function can never be convex,
the assertion follows.
\end{proof}

\begin{proposition}
\label{polprops}
Let $\kegela$ be a pointed cone in 
$\gittera_{\Q}$.
Then we have the following statements 
for the associated groups of 
$\kegela$-polyhedra: 
\begin{enumerate}
\item 
There is a commutative diagram of 
canonical, injective homomorphisms
of monoids:
$$
\xymatrix{
\Pol^+_{\kegela}(\gittera) \ar[d] \ar[r]
&
\Pol^+_{\kegela}(\gittera_{\Q})  \ar[d]
\\
\Pol_{\kegela}(\gittera) \ar[r]
&
\Pol_{\kegela}(\gittera_{\Q}) 
}
$$
\item
The multiplication of elements
$\polyedera \in \Pol^+_{\kegela}(\gittera_{\Q})$
by positive rational numbers
$\alpha \in \Q_{> 0}$, defined as 
\begin{eqnarray*} 
\alpha \polyedera 
& := &
\{\alpha v; \; v \in \polyedera\},
\end{eqnarray*}
uniquely
extends to a scalar multiplication 
$\Q \times \Pol_{\kegela}(\gittera_{\Q}) 
\to \Pol_{\kegela}(\gittera_{\Q})$ 
making $\Pol_{\kegela}(\gittera_{\Q})$
into a rational vector space.
\item
The group $\Pol_{\kegela}(\gittera)$ of 
integral $\kegela$-polyhedra is 
free abelian, and we have a canonical
isomorphism
\begin{eqnarray*} 
\Pol_{\kegela}(\gittera_{\Q})
& \cong & 
\Q \otimes_{\Z} \Pol_{\kegela}(\gittera).
\end{eqnarray*}
\item
For every element $u \in \kegela^\vee$,
there is a unique 
linear evaluation functional 
$\eval_{u}  \colon  
\Pol_{\kegela}(\gittera_{\Q}) 
\to  \Q$
satisfying 
\begin{eqnarray*}
\eval_{u}(\polyedera)
& = & 
\min_{v \in \polyedera} \langle u, v \rangle,
\quad \text{if }
\polyedera \in \Pol^+_{\kegela}(\gittera_{\Q}).
\end{eqnarray*}
\item 
Two elements 
$\polyedera_1, \polyedera_2 \in \Pol_{\kegela}(\gittera_{\Q})$
coincide if and only if 
$\eval_{u}(\polyedera_1) = \eval_{u}(\polyedera_2)$
holds for all  $u \in \kegela^\vee$.
\item
An element $\polyedera \in \Pol_{\kegela}(\gittera_{\Q})$
belongs to $\Pol_{\kegela}(\gittera)$ if and only
if for every $u \in \kegela^\vee \cap \dualgittera$,
the evaluation $\eval_u(\polyedera)$ is an integer.
\end{enumerate}
\end{proposition}

\begin{proof}
For assertion~(i), note that by Proposition~\ref{cancel}, 
the monoids
$\Pol^+_{\kegela}(\gittera)$ and
$\Pol^+_{\kegela}(\gittera_{\Q})$
embed into their Gro\-then\-dieck groups.
The rest of the assertion is
a consequence of functoriality of
the Grothendieck group.
Similarly, existence and uniqueness 
of the scalar multiplication in 
assertion~(ii) can be established
via functoriality of the Grothen\-dieck
group.

For assertion~(iii), note that the map 
$\Pol^+_{\kegela}(\gittera_{\Q}) 
\to 
\CPL_{\Q}(\kegela^{\vee})$ 
of Proposition~\ref{poly2suppfunct} sends 
the elements of $\Pol^+_{\kegela}(\gittera)$
to functions having integer values on
$\kegela^{\vee} \cap \dualgittera$.
Thus, we may view $\Pol_{\kegela}(\gittera)$
as a subgroup of the abelian group of 
all integer-valued functions on the countable 
set $\kegela^{\vee} \cap \dualgittera$.
Countable subgroups of this group are 
free abelian, see for example~\cite[Satz~1]{Speck}.
This applies to $\Pol_{\kegela}(\gittera)$.
The claimed isomorphism is then easily 
obtained by considering a $\Z$-basis for
$\Pol_{\kegela}(\gittera)$.

On $\Pol^+_{\kegela}(\gittera_\Q)$,
the existence of the evaluation functional
asserted in~(iv) is due to 
Proposition~\ref{poly2suppfunct}; in fact, we have 
$\eval_u(\polyedera) = h_\polyedera(u)$.
The unique continuation to $\Pol_{\kegela}(\gittera_\Q)$,
is, once more, a consequence of the universal
property of the Grothendieck group.

To verify the ``if'' part of assertion~(v), write 
$\polyedera_i = \polyedera_i^+ -  \polyedera_i^-$
with two $\kegela$-polyhedra 
$\polyedera_i^+$ and $\polyedera_i^-$.
Then the sums
$
\polyedera_1 + \polyedera_1^- +  \polyedera_2^-
$
as well as 
$
\polyedera_2 + \polyedera_1^- +  \polyedera_2^-
$
are $\kegela$-polyhedra,
and all their evaluations coincide.
Thus, Proposition~\ref{poly2suppfunct}
says that these two $\kegela$-polyhedra
coincide. The assertion follows.

For the ``if'' part of assertion~(vi), 
it suffices to
consider $\kegela$-polyhedra $\polyedera$, 
because any element of 
$\Pol_{\kegela}(\gittera_\Q)$
can be shifted into 
$\Pol^+_{\kegela}(\gittera_\Q)$
by adding an integral $\kegela$-polyhedron.
For any vertex $v \in \polyedera$, the linear
forms $u \in \kegela^\vee  \cap \dualgittera$ 
attaining their minimum over $\polyedera$
in $v$ generate $\dualgittera$ as a lattice,
because the cone of $\gkfaecher(\polyedera)$ 
corresponding to $v$ is of full dimension.
Hence, the vertices of $\polyedera$ belong to 
$\gittera$ if all evaluations $\eval_u$,
where $u \in \kegela^\vee \cap \dualgittera$,
are integral on $\polyedera$. 
\end{proof}

\section{Polyhedral divisors}\label{section1}

In this section, we introduce the language of 
polyhedral divisors.
The idea is to allow not only integral or 
rational numbers as coefficients of a divisor, 
but more generally, integral or arbitrary
tailed polyhedra.
The essential points of this section are the
definition of proper polyhedral divisors
(pp-divisors)
and an interpretation of this notion in terms
of convex piecewise linear maps, see~\ref{ppoldivdef}
and~\ref{cpl}.

Here, and moreover in the entire paper, 
the words algebraic variety refer to an integral 
scheme of finite type over a variety over 
an algebraically closed field $\K$ of 
characteristic zero (though we expect to hold
the results as well in positive characteristics,
with basically the same proofs). 
By a point, we always mean a closed point,
and $\K(Y)$ denotes the function field of $Y$.

The following class of varieties
will be of special importance for us; it
comprises the affine as well as the 
projective ones, 
compare also~\cite{hyperkaehler}.

\begin{definition}
An algebraic variety $Y$ is said to be 
{\em semiprojective\/} if its $\K$-algebra
of global functions
$A_{0} := \Gamma(Y,\mathcal{O})$ is finitely generated, 
and $Y$ is projective over $Y_{0} := \Spec(A_{0})$.
\end{definition}

The groups of Weil and Cartier divisors 
on a {\em normal\/} algebraic variety $Y$ are
denoted by $\WDiv(Y)$ and $\CDiv(Y)$,
and the corresponding vector spaces 
of rational divisors are denoted by 
$\WDiv_\Q(Y)$ and $\CDiv_\Q(Y)$.
Since $Y$ is normal, we have the inclusions
$\CDiv(Y) \subset \WDiv(Y)$, and 
$\CDiv_\Q(Y) \subset \WDiv_\Q(Y)$.

Let us briefly recall the basic notions around
divisors used later. 
The sheaf of sections 
$\mathcal{O}(D)$ of a rational 
Weil divisor $D$ on 
a normal algebraic variety $Y$ is,
similar to the usual case,
defined via
$$
\begin{array}{ccccc}
\Gamma(V,\mathcal{O}(D))
& := &
\{f\in \K(Y) \kst \ddiv(f\vert_{V}) + D\vert_{V} \geq 0\}
& = & 
\Gamma(V,\mathcal{O}(\rd{D})), 
\end{array}
$$
where $V \subset Y$ is open and $\rd{D}$ denotes the 
round-down divisor of $D$.
For a section 
$f \in \Gamma(Y,\mathcal{O}(D)) \subset \K(Y)$ 
of a rational Weil divisor $D$ on a normal 
algebraic variety $Y$, we define its 
{\em zero set\/} and its
{\em non-vanishing locus\/} as
$$
Z(f) \; := \; \Supp(\ddiv(f) + D),
\qquad 
Y_{f} \; := \; Y \setminus Z(f).
$$
Moreover, $D \in \CDiv_\Q(Y)$ is called 
{\em semiample\/}
if it admits a basepoint-free multiple, i.e.
for some $n \in \Z_{> 0}$ 
the sets $Y_f$, where 
$f \in \Gamma(Y,\mathcal{O}(nD))$,
cover~$Y$.
We also need a straightforward generalization
of the concept of a big divisor on a projective 
variety, compare~\cite[Lemma.~2.60]{KoMo}.

\begin{definition}
We say that a divisor $D \in \CDiv_\Q(Y)$ on 
a variety $Y$ is {\em big\/} if for some 
$n \in \Z_{> 0}$  there is
a section $f \in \Gamma(Y,\mathcal{O}(nD))$
with an affine non-vanishing locus $Y_f$. 
\end{definition}

Now we turn to divisors with tailed polyhedra as coefficients.
Here are the first definitions.

\begin{definition}
\label{poldivdef}
Let $Y$ be a normal algebraic variety; 
let $\gittera$ be a lattice, 
and let $\kegela \subset \gittera_{\Q}$ 
be a pointed cone.
\begin{enumerate}
\item
The groups of {\em rational polyhedral Weil divisors\/}
and {\em rational polyhedral Cartier divisors\/} 
of $Y$ with respect to $\kegela \subset \gittera_{\Q}$ 
are
\begin{eqnarray*}
\WDiv_\Q(Y,\kegela)  
& :=  &
\Pol_\kegela(\gittera_{\Q}) \otimes_\Z \WDiv(Y), 
\\
\CDiv_\Q(Y,\kegela)  
& :=  &
\Pol_\kegela(\gittera_{\Q}) \otimes_\Z \CDiv(Y).
\end{eqnarray*}
\item
The groups of {\em integral polyhedral Weil divisors\/}
and {\em integral polyhedral Cartier divisors\/} 
of $Y$ with respect to $\kegela \subset \gittera_{\Q}$ 
are
\begin{eqnarray*}
\WDiv(Y,\kegela)  
& :=  &
\Pol_\kegela(\gittera) \otimes_\Z \WDiv(Y), 
\\
\CDiv(Y,\kegela)  
& :=  &
\Pol_\kegela(\gittera) \otimes_\Z \CDiv(Y).
\end{eqnarray*}
\end{enumerate}
\end{definition}

>From here on, when we speak of 
divisors, or polyhedral divisors
we mean rational ones; 
if we want to consider integral divisors,
then this is explicitly stated.
Here is a list of first properties
of the groups of polyhedral divisors.

\begin{proposition}
\label{poldiveval}
Let $Y$ be a normal algebraic variety; 
let $\gittera$ be a lattice, 
and let $\kegela \subset \gittera_{\Q}$ 
be a pointed cone.
\begin{enumerate}
\item
$\WDiv_\Q(Y,\kegela)$ and $\CDiv_\Q(Y,\kegela)$ 
are rational vector spaces,
and $\WDiv(Y,\kegela)$ and $\CDiv(Y,\kegela)$ 
are free abelian groups.
\item 
There is commutative diagram of canonical 
injections
$$
\xymatrix{
\CDiv(Y,\kegela) \ar[r]\ar[d]
&
\WDiv(Y,\kegela) \ar[d] \\
\CDiv_\Q(Y,\kegela) \ar[r]
&
\WDiv_\Q(Y,\kegela) \;\; .
}
$$
Moreover, we have canonical isomorphisms
\begin{eqnarray*}
\WDiv_\Q(Y,\kegela)
& \cong &
\Q \otimes_\Z 
\WDiv(Y,\kegela), 
\\
\CDiv_\Q(Y,\kegela)
& \cong &
\Q \otimes_\Z 
\CDiv(Y,\kegela).
\end{eqnarray*}
\item
For every element $u \in \kegela^\vee$,
there is a well defined linear evaluation functional 
\begin{eqnarray*}
\WDiv_\Q(Y,\kegela) 
& \to & 
\WDiv_\Q(Y),
\\
\poldiva \, = \; \sum \polyedera_{i} \otimes D_{i} 
& \mapsto &
\poldiva(u) \; := \; \sum \eval_{u}(\polyedera_{i}) D_{i}.
\end{eqnarray*}
\item
Two polyhedral divisors 
$\poldiva_1, \poldiva_2 \in \WDiv_\Q(Y,\kegela)$
coincide if and only if we have  
$\poldiva_1(u) = \poldiva_2(u)$ for 
all 
$u \in \kegela^\vee$.
\item
A polyhedral divisor
$\poldiva \in \WDiv_\Q(Y,\kegela)$
is integral if and only if 
all its evaluations 
$\poldiva(u)$, where 
$u \in \kegela^\vee \cap \dualgittera$,
are integral divisors.
\item
A polyhedral divisor
$\poldiva \in \WDiv_\Q(Y,\kegela)$
is Cartier if and only if 
all its evaluations 
$\poldiva(u)$, where 
$u \in \kegela^\vee$,
are Cartier divisors.
\end{enumerate}
\end{proposition}

\begin{proof}
The first three assertions are immediate
consequences of Proposition~\ref{polprops}.
For the last three assertions, note first 
that any $\poldiva \in \WDiv_\Q(Y,\kegela)$
allows a representation 
$$
\begin{array}{ccccc}
\poldiva 
& = & 
\sum \polyedera_{i} \otimes D_{i}
& = & 
\sum (\polyedera_{i}^+-\polyedera_{i}^-) \otimes D_{i}
\end{array}
$$
with prime divisors $D_i$, 
and coefficients
$\polyedera_i \in \Pol(\gittera_\Q)$
which are of the form 
$\polyedera_i = \polyedera_i^+ - \polyedera_i^-$
with  $\kegela$-polyhedra
$\polyedera_i^+, \polyedera_i^-
\in \Pol_\kegela^+(\gittera_\Q)$. 

Combining this observation with 
Proposition~\ref{polprops}~(v) gives
assertion~(iv).
Moreover, if all evaluations 
$\poldiva(u)$, where $u \in \kegela^\vee \cap \dualgittera$,
are integral divisors, then 
Proposition~\ref{polprops}~(vi) yields
that the coefficients $\polyedera_i$ 
are integral.
This merely means that $\poldiva$ is an integral
polyhedral divisor.

Now, suppose that all evaluations of $\poldiva$
are Cartier divisors.
Consider the vector space
$W \subset \WDiv_\Q(Y)$ generated 
by the prime divisors $D_i$ and 
the vector subspace $W_0 \subset W$ 
generated by the evaluations 
$\poldiva(u)$, where 
$u \in \kegela^\vee \cap \dualgittera$.
Then, with respect to a basis $E_1, \ldots, E_k$ 
of $W_0$, we may write
\begin{eqnarray*}
\poldiva(u)
& = &
\sum g_j^+(u)E_j - g_j^-(u)E_j,
\end{eqnarray*}
where the $g_j^+$ and the $g_j^-$ 
are nonnegative linear combinations
of the support functions corresponding
to $\polyedera_i^+$ and $\polyedera_i^-$.
In particular, $g_j^+$ and $g_j^-$ are 
piecewise linear, convex functions. 
Thus, by assertion~(iv),
we have a representation
\begin{eqnarray*}
\poldiva
& = & 
\sum (\polyederb_j^+ - \polyederb_j^-) \otimes E_j,
\end{eqnarray*}
where $\polyederb_j^+$ and $\polyederb_j^-$ are the 
$\kegela$-polyhedra having the functions 
$g_j^+$ and $g_j^-$ as their support functions,
respectively. Since the $E_j$ are 
Cartier divisors, this shows that $\poldiva$
is a polyhedral Cartier divisor.
\end{proof}

In the spirit of the last two statements of this
proposition, we may introduce further notions
for polyhedral divisors.

\begin{definition}
Let $Y$ be a normal algebraic variety; 
let $\gittera$ be a lattice, 
and let $\kegela \subset \gittera_{\Q}$ 
be a pointed cone.
\begin{enumerate}
\item
We call a polyhedral divisor $\poldiva \in \WDiv_\Q(Y,\kegela)$  
{\em effective (written $\poldiva \ge 0$)\/}, if all evaluations 
$\poldiva(u)$, where $u \in \kegela^\vee$, are effective divisors.
\item 
We call a polyhedral divisor 
$\poldiva \in \CDiv_\Q(Y,\kegela)$
{\em semiample\/} if all evaluations 
$\poldiva(u)$, where $u \in \kegela^\vee$, 
are semiample divisors.
\end{enumerate}
\end{definition}

\begin{example}
\label{usual2polyhedral}
Let $Y$ be any normal variety, $\gittera := \Z$, 
and $\kegela := \Q_{\ge 0}$.
Then we have a canonical isomorphism
$$
\Div_{\Q}(Y)
\; \to \; 
\Div_{\Q}(Y,\sigma),
\qquad
\sum \alpha_{i} D_{i}
\; \mapsto \; 
\sum (\alpha_{i} + \sigma) \otimes D_{i} .
$$
Integral (effective, Cartier, semiample) 
divisors correspond to
integral (effective, Cartier, semiample) 
polyhedral divisors.
The inverse isomorphism is
given by
$$
\Div_{\Q}(Y,\sigma)
\; \to \; 
\Div_{\Q}(Y),
\qquad
\poldiva \; \mapsto \; \poldiva(1).
$$
\end{example}

We come to the central definition of the 
paper; we introduce the class 
of proper polyhedral divisors.

\begin{definition}
\label{ppoldivdef}
Let $Y$ be a normal algebraic variety; 
let $\gittera$ be a lattice, 
and let $\kegela \subset \gittera_{\Q}$ 
be a pointed cone.
A {\em proper polyhedral divisor\/}
(abbreviated {\em pp-divisor\/})
on $Y$ with respect to $\kegela \subset N_\Q$
is a semiample polyhedral divisor 
$\poldiva \in \CDiv_{\Q}(Y,\kegela)$
such that
\begin{enumerate}
\item 
there is a representation
$\poldiva = \sum \polyedera_i \otimes D_i$ with
$\kegela$-polyhedra $\polyedera_i \in \Pol^{+}_{\kegela}(\gittera_\Q)$
and effective divisors $D_i \in \WDiv(Y)$,
\item 
for every $u \in \relint{\kegela^\vee}$,
the evaluation $\poldiva(u) \in \CDiv_{\Q}(Y)$ 
is a big divisor on $Y$.
\end{enumerate}
\end{definition}

Clearly, the sum of two pp-divisors with respect to 
a given $\kegela \subset \gittera_\Q$ is again
a pp-divisor with respect to $\kegela \subset \gittera_\Q$ . 
Thus, these polyhedral divisors form a semigroup.
Our notation is the following.

\begin{definition}
Let $Y$ be a normal algebraic variety; 
let $\gittera$ be a lattice, 
and let $\kegela \subset \gittera_{\Q}$ 
be a pointed cone. 
The semigroup of all pp-divisors on $Y$
with respect to $\kegela \subset \gittera_\Q$ 
is denoted by $\PPDiv_{\Q}(Y,\kegela)$.
\end{definition}

We show now that the pp-divisors 
$\poldiva \in \PPDiv_{\Q}(Y,\kegela)$ correspond 
to certain convex piecewise linear
maps $\kegela^{\vee} \to \CDiv_{\Q}(Y)$.
The precise definition of these maps is the following.

\begin{definition}
\label{cpldef}
Let $Y$ be a normal variety; let $\dualgittera$ 
be a lattice, and let
$\dualkegela \subset \dualgittera_\Q$ be a cone 
of full dimension.
We say that a map 
$\cpla \colon \dualkegela \to \CDiv_{\Q}(Y)$ is
\begin{enumerate}
\item \label{cpldef1}
 {\em convex\/} if $\cpla(u)+\cpla(u')\le\cpla(u+u')$
  holds for any two elements $u,u' \in \dualkegela$,
\item \label{cpldef2}
 {\em piecewise linear\/} if there is a quasifan $\gkfaecher$
 in $\dualgittera_\Q$ having $\dualkegela$ as its support such that 
  $\cpla$ is linear on the cones of $\gkfaecher$,
\item \label{cpldef3}
  {\em strictly semiample\/} if $\cpla(u)$ is always semiample
  and, for $u \in \relint{\dualkegela}$, it is even big.
\end{enumerate}
\end{definition}

The sum of two convex, piecewise linear, 
strictly semiample maps is again of this type. 
We use the following notation:

\begin{definition}
Let $Y$ be a normal variety; let $\dualgittera$ 
be a lattice, and let
$\dualkegela \subset \dualgittera_\Q$ be a cone 
of full dimension.
The semigroup of all convex, piecewise linear, strictly semiample 
maps $\cpla \colon \dualkegela \to \CDiv_{\Q}(Y)$ is denoted
by $\CPL_{\Q}(Y,\dualkegela)$.
\end{definition}

With these definitions, we are ready to state 
the analogue of Proposition~\ref{poly2suppfunct}
for proper polyhedral divisors.

\begin{proposition}
\label{cpl}
Let $Y$ be a normal variety, $\gittera$
a lattice, and 
$\kegela \subset \gittera_\Q$ a pointed cone.
Then there is a canonical isomorphism of 
semigroups:
$$
\PPDiv_{\Q}(Y,\kegela)
\; \to \; 
\CPL_{\Q}(Y,\kegela^\vee),
\qquad
\poldiva \; \mapsto \; \left[u \mapsto \poldiva(u) \right].
$$
Under this isomorphism, the integral polyhedral divisors
correspond to those maps sending 
$\kegela^{\vee} \cap M$
to $\CDiv(Y)$.
\end{proposition}

\begin{proof}
By Proposition~\ref{poldiveval}, 
the assignment is a well defined
injective homomorphism. 
Thus, we only need to verify that any convex
piecewise linear map 
$\cpla \colon \kegela^\vee \to \CDiv_\Q(Y)$,
in the sense of
Definition~\ref{cpldef},
arises from a polyhedral divisor.
Since $\kegela^\vee$ is polyhedral,
there occur only finitely many prime 
divisors $D_1, \ldots, D_r$ in the images
$\cpla(u)$, where $u \in \sigma^\vee$.
Thus, we may write
\begin{eqnarray*}
\cpla(u) 
& = &
\sum_{i=1}^{r} h_i(u) D_i,
\end{eqnarray*}
where every $h_i \colon \sigma^\vee \to \Q$ is a
$\Q$-valued, convex, piecewise linear 
function in the usual sense.
According to Proposition~\ref{poly2suppfunct},
each of the functions $h_i$ is of the form 
$u \mapsto \eval_u(\polyedera_i)$ with a 
$\kegela$-polyhedron $\polyedera_i \subset N_\Q$.
Consequently, the sum of all $\polyedera_i \otimes D_i$ 
is a polyhedral divisor 
defining the map~$\cpla \colon \kegela^\vee \to \CDiv_\Q(Y)$.
\end{proof}

This observation allows us to switch freely between
pp-divisors and convex, piecewise linear,
strictly semiample maps. In particular, we 
denote these objects by the same symbol, preferably
by the gothic letter $\poldiva$.

\begin{example}
\label{exsa}
Let $Y$ be a smooth projective curve;
let $\gittera$ be a lattice, and let 
$\kegela \subset \gittera_\Q$ be a
pointed cone.
To any polyhedral divisor on $Y$ 
with respect to $\kegela \subset \gittera_\Q$,
we associate its {\em polyhedral degree\/} by setting
$$
\begin{array}{ccccc}
\deg\left(\sum  \polyedera_i \otimes D_i \right)
& := &
\sum \deg(D_i) \polyedera_i 
& \in & 
\Pol_{\kegela}(\gittera_\Q).
\end{array}
$$
This does not depend on the representation
of a given $\poldiva \in \CDiv_{\Q}(Y,\kegela)$,
and $\eval_u(\deg(\poldiva))$
equals $\deg(\poldiva(u))$
for any $u \in \kegela^\vee$.
We will figure out, in terms of the degree,
when a given $\poldiva \in \CDiv_{\Q}(Y,\kegela)$
is a pp-divisor. 

First, recall that on the curve $Y$, 
a divisor is big if and only
if it has positive degree, and a divisor is 
semiample if and only if it is big or
some multiple of it is a principal divisor. 
Consequently, $\poldiva \in \CDiv_{\Q}(Y,\kegela)$
is a pp-divisor if and only if the following
holds:
\begin{enumerate}
\item 
$\poldiva = \sum \polyedera_i \otimes \{y_i\}$,
with $y_i \in Y$ pairwise disjoint
and $\polyedera_i \in \Pol_\kegela^+(N_\Q)$,
\item
the $\kegela$-polyhedron $\deg(\poldiva)$ 
is a proper subset of the cone $\kegela$,
\item 
if $\eval_u(\deg(\poldiva)) = 0$, 
then 
$u\in\partial\kegela^\vee$ and
a multiple of $\poldiva(u)$ is principal.
\end{enumerate} 
Note that the first of these conditions is a 
reformulation of Condition~\ref{ppoldivdef}~(i).
Moreover, the last two conditions are  
satisfied if $\deg(\poldiva)$ is contained
in the relative interior of $\kegela$.
\end{example}

\section{pp-divisors and torus actions}
\label{section2}

In this section, we formulate the first results 
of this paper.
They show that the affine normal varieties with 
an effective algebraic torus action arise from 
pp-divisors on normal semiprojective 
varieties.

Let us briefly fix our notation around torus 
actions and also recall a little background.
An (algebraic) torus is an affine algebraic group
$\torusa = \Spec(\K[\dualgittera])$,
where $\dualgittera$ is a lattice, and 
$\K[M]$ denotes the associated group algebra.
For an element $u \in \dualgittera$, 
we denote, as usual, the corresponding 
character by 
$\chi^u \colon \torusa \to \K^*$.

If a torus $\torusa = \Spec(\K[\dualgittera])$ 
acts on a variety $X$, 
then we always assume that this action is given 
by a morphism 
$$
\torusa \times X \; \to \; X,
\qquad
(t,x) \mapsto t \mal x,
$$
and we also speak about the $T$-variety $X$.
A {\em semiinvariant\/} with respect
to the character $\chi^u \colon \torusa \to \K^*$
is a function $f \in \Gamma(X,\mathcal{O})$
satisfying
$$ 
f(t \mal x)
\;  = \;  
\chi^{u}(t) f(x)
\quad
\text{for all }
(t,x) \in T \times X.
$$

We write $\Gamma(X,\mathcal{O})_u$ 
for the vector space of semiinvariants 
with respect to $\chi^u$, 
and $\Gamma(X,\mathcal{O})^T$ 
for the algebra of invariants, 
i.e.\ the semiinvariants with 
respect to $\chi^0$.
The action of $\torusa$ on $X$ is called effective
if the neutral element of $\torusa$ is the only 
element acting trivially on $X$.

A morphism $\pi \colon X \to Y$ is called a
{\em good quotient\/} for a $\torusa$-action
on $X$ if it is affine, $\torusa$-invariant, 
i.e.\ constant on $T$-orbits, and the pullback map   
$\pi^* \colon \mathcal{O}_Y \to \pi_*(\mathcal{O}_X)^\torusa$
is an isomorphism.
If a good quotient exists, then it is unique
up to an isomorphism, 
and the quotient space is frequently denoted as
$X \quot T$.

The possible actions of a torus 
$\torusa  = \Spec(\K[\dualgittera])$ 
on an affine variety
$X = \Spec(A)$ correspond to $\dualgittera$-gradings
of the algebra $A$: given a $\torusa$-action, 
the homogeneous part $A_u \subset A$ for $u \in M$
consists precisely of the semiinvariants with respect
to the character $\chi^u \colon \torusa \to \K^*$.

The {\em weight monoid\/} $S$ of a $\torusa$-action 
on $X = \Spec(A)$ consists of all 
$u \in \dualgittera$ with $A_u \ne \{0\}$,
and the {\em weight cone\/} is the 
(convex, polyhedral) cone
$\gewkegel \subset \dualgittera_{\Q}$ 
generated by $S$.
We will usually denote the $\dualgittera$-grading 
of $A$ as 
\begin{eqnarray*}
A & = & \bigoplus_{u \in \gewkegel \cap \dualgittera} A_{u}. 
\end{eqnarray*}

Let us present the first result.
Fix a normal semiprojective variety~$Y$, 
a lattice $\gittera$, a pointed cone 
$\kegela \subset \gittera_{\Q}$, and a pp-divisor 
$\poldiva \in \PPDiv_{\Q}(Y,\kegela)$.
Then these data define a sheaf of multigraded 
algebras on $Y$:
the convexity property~\ref{cpldef}~(i)
of the map $u \mapsto \poldiva(u)$
ensures the existence of canonical multiplication 
maps
\begin{eqnarray*}
 \mathcal{O}(\poldiva(u)) 
   \otimes \mathcal{O}(\poldiva(u')) 
 & \to &
\mathcal{O}(\poldiva(u+u')),
\end{eqnarray*}
and thus, the sheaves $\mathcal{O}(\poldiva(u))$,
where $u \in \kegela^\vee \cap \dualgittera$,
can be put together to an 
$\mathcal{O}_{Y}$-algebra $\mathcal{A}$, 
graded by the monoid $\kegela^\vee \cap \dualgittera$.
Now we take the relative spectrum
$\til{X} := \Spec_Y(\mathcal{A})$.
Here are the basic properties of this
construction.

\begin{theorem}
\label{data2variety}
Let $Y$ be a normal semiprojective variety, 
$\gittera$ a lattice, 
$\kegela \subset \gittera_{\Q}$ 
a pointed cone,
and $\dualgittera := \Hom(\gittera,\Z)$.
Given a pp-divisor
$\poldiva \in \PPDiv_{\Q}(Y,\kegela)$, 
consider the $\mathcal{O}_{Y}$-algebra
$$ 
\begin{array}{ccccc}
\mathcal{A}
& := &
\displaystyle
\bigoplus_{u \in \kegela^{\vee} \cap \dualgittera} \mathcal{A}_{u}
& := & 
\displaystyle
\bigoplus_{u \in \kegela^{\vee} \cap \dualgittera} 
             \mathcal{O}(\poldiva(u)),
\end{array}
$$
the algebraic torus $\torusa := \Spec(\K[\dualgittera])$,
and the relative spectrum
$\til{X} = \Spec_Y(\mathcal{A})$.
Then the following statements hold.
\begin{enumerate}
\item The scheme $\til{X}$ is a normal
   algebraic variety of dimension $\dim(Y)+\dim(\torusa)$, 
   and the grading of $\mathcal{A}$ defines an effective torus 
   action $\torusa \times \til{X} \to \til{X}$
   having the canonical map $\quotabb \colon \til{X} \to Y$ as a 
   good quotient.
\item The ring of global sections 
   $A := 
   \Gamma(\til{X},\mathcal{O})=
   \Gamma(Y,\mathcal{A}) 
   $ 
   is a finitely generated $M$-graded normal $\K$-algebra, 
   and we have a proper, birational
   $\torusa$-equivariant contraction morphism $\til{X} \to X$
   with $X := \Spec(A)$.
\item
   Let $u \in \kegela^{\vee} \cap M$ and $f \in A_{u}$. Then we have
   $\pi(\til{X}_{f}) = Y_{f}$. In particular, if $Y_{f}$
   is affine, then so is $\til{X}_{f}$, and the canonical map
   $\til{X}_f\to X_f$ is an isomorphism.
   Moreover, even for non-affine $Y_f$, we have
  \begin{eqnarray*}
   \Gamma(Y_f,\mathcal{A}) 
     &  = & 
   \bigoplus_{u\in \kegela^{\vee} \cap \dualgittera} (A_f)_u.
   \end{eqnarray*}
\end{enumerate}%
\end{theorem}

The proof will be given in Section~\ref{firstproof}.
As the following two examples show,
the result extends both the construction of affine toric 
varieties, see e.g.~\cite{Fulton},
and the Dolgachev-Demazure-Pinkham 
construction of good $\K^{*}$-actions,
see~\cite{Dol}, \cite{Dem} and~\cite{Pi}:

\begin{example}[Affine Toric Varieties]
\label{afftorvar}
Let $N$ be a lattice and $\kegela \subset \gittera_\Q$
a cone.
Then the associated affine toric variety $X_{\kegela}$
is defined as the spectrum of the 
semigroup algebra $\K[\kegela^{\vee} \cap \dualgittera]$
where $\dualgittera$ is the dual lattice of $\gittera$, 
and $\kegela^{\vee}$ is the dual cone of~$\kegela$.
We can also obtain $X_\kegela$ as the $\til{X}=X$ 
of a pp-divisor:
Take $Y := \Spec(\K)$, 
and let
$\poldiva \in \PPDiv_{\Q}(Y,\kegela)$ 
be the trivial divisor.
\end{example}

\begin{example}[Good $\K^{*}$-Actions]
\label{good}
Let $Y$ be a normal projective variety,
and let $D$ be an ample rational
Cartier divisor on $Y$.
These data give rise to a pp-divisor: 
take $N := \Z$; let 
$\kegela \subset \Q$ be the positive ray, 
and consider 
$\poldiva := (1 + \kegela) \otimes D$.
Then $\poldiva$ corresponds to the map
$$
\kegela^{\vee} \; \to \; \CDiv_{\Q}(Y),
\qquad
u \; \mapsto \; uD.
$$

If $D$ is an integral Cartier divisor,
then $\til{X} = \Spec_Y(\mathcal{A})$ 
is the total space of a line bundle, and $\til{X} \to X$
is the $\K^{*}$-equivariant contraction of the zero
section. 
Thus, the affine variety $X$ is an affine cone 
over $Y$. For $D$ being a rational divisor, 
$X$ is usually called a generalized cone over $Y$.
\end{example}

In our second result, we go the other way around.
We show that every normal affine 
variety $X$ with an effective torus action arises 
from a pp-divisor 
$\poldiva \in \PPDiv_\Q(Y,\kegela)$ 
on some normal semiprojective variety $Y$
in the sense of the preceding construction.

\begin{theorem}\label{variety2data}%
Let $X$ be a normal affine variety and suppose
that  $\torusa = \Spec(\K[\dualgittera])$
acts effectively on $X$ with weight cone 
$\dualkegela \subset \dualgittera_{\Q}$.
Then there exists a normal semiprojective
variety $Y$ and a pp-divisor
$\poldiva \in \PPDiv_{\Q}(Y,\dualkegela^{\vee})$
such that we have an isomorphism of graded 
algebras:
\begin{eqnarray*}
\Gamma(X,\mathcal{O})
&
\cong
&
\bigoplus_{u \in  \gewkegel  \cap \dualgittera} 
\Gamma(Y,\mathcal{O}(\poldiva(u))).
\end{eqnarray*}
\end{theorem}

For the proof, we refer to Section~\ref{limit}.
Our construction of the semiprojective variety~$Y$
and the pp-divisor $\poldiva$ 
is basically canonical. It relies on the chamber structure 
of the collection of all GIT-quotients of $X$ 
that arise from possible linearizations of 
the trivial bundle.

We conclude this section with a further example.
We indicate how to recover the Flenner--Zaidenberg 
description~\cite{FlZa} of normal affine $\K^{*}$-surfaces 
from Theorems~\ref{data2variety} and~\ref{variety2data}:

\begin{example}[Normal affine $\K^{*}$-surfaces]
\label{XXX}
Any normal affine surface $X$ with effective 
$\K^{*}$-action arises from a pp-divisor
$\poldiva \in \PPDiv_{\Q}(Y,\sigma)$ with 
a normal and hence smooth curve~$Y$ and a 
pointed cone $\sigma \subset \Q$ where
the lattice is $\Z$. 

The curve $Y$ is either affine or projective.
In the latter case, we have $\sigma \ne \{0\}$,
because otherwise the convexity property
$\poldiva(1) + \poldiva(-1) \le \poldiva(0)=0$ would 
contradict strict semiampleness of $\poldiva$.
Thus, up to switching the action by
$t \mapsto t^{-1}$, there are three cases:

\begin{ecase}
The curve $Y$ is projective and  $\sigma = \Q_{\ge0}$ holds. 
Then the $\K^{*}$-action on 
$X$ is good, i.e.\ it has an (isolated) attractive fixed point. 
Moreover, $\poldiva$ is of the form
\begin{eqnarray*}
\poldiva 
& = &
\sum_{i=1}^{r} [v_{i},\infty) \otimes \{y_i\}, 
\end{eqnarray*} 
with $y_{i} \in Y$ and $v_{i} \in \Q$ such that
$v_{1} + \ldots + v_{r} > 0$.
Note that in this case, $\poldiva$ is determined by its 
evaluation $D := \poldiva(1)$,
namely $\poldiva=[1,\infty)\otimes D$.
\end{ecase}
\begin{pcase}
The curve $Y$ is affine, and $\sigma = \Q_{\ge0}$ holds. 
Then the $\K^{*}$-action has an attractive fixed
point curve, isomorphic to $Y$. Moreover,
\begin{eqnarray*} 
\poldiva 
& = &
\sum_{i=1}^{r} [v_{i},\infty)  \otimes \{y_i\},
\end{eqnarray*} 
with $y_{i} \in Y$, but no condition on the numbers
$v_{1}, \ldots, v_{r} \in \Q$.
Again, $\poldiva$ is determined by its evaluation 
$D := \poldiva(1)$.
\end{pcase}
\begin{hcase}
The curve $Y$ is affine, and $\sigma = \{0\}$ holds. 
Then the generic $\K^{*}$-orbit is closed. For the pp-divisor 
$\poldiva$, we obtain a representation
\begin{eqnarray*} 
\poldiva 
& = & 
\sum_{i=1}^{r} [v_{i}^{-},v_{i}^{+}]  \otimes \{y_i\}, 
\end{eqnarray*} 
with $y_{i} \in Y$, and $v_{i}^{-} \le v_{i}^{+}$. 
Note that $\poldiva$ is determined by 
$D^{-} := \poldiva(-1)$ and $D^{+} := \poldiva(1)$.
This pair satisfies $D^{-} + D^{+} \le 0$, 
and we have 
\begin{eqnarray*}
\poldiva
& = & 
\{1\} \otimes D^+
\; - \; 
[0,1] \otimes (D^- + D^+).
\end{eqnarray*}
%
\end{hcase}
\end{example}

\section{Proof of Theorem~\ref{data2variety}}
\label{firstproof}

This section is devoted to proving 
Theorem~\ref{data2variety}.
Note that parts of the assertions~(i) and~(ii) 
are well known for the case that $Y_0$ is a point 
and $u \mapsto \poldiva(u)$ is linear, see for 
example~\cite[Thm.~4.2]{Zar} and~\cite[Lemma~2.8]{HuKe}. 

We start with a basic observation concerning 
multigraded rings which will also be used apart
from the proof of Theorem~\ref{data2variety}. 

\begin{lemma}\label{klaustrick}
Let $\dualgittera$ be a lattice, and let $A$ be a 
finitely generated $\dualgittera$-graded $\K$-algebra. 
Then, every (convex, polyhedral) cone 
$\dualkegela \subset \dualgittera_{\Q}$ 
defines a finitely generated $\K$-algebra
\begin{eqnarray*} 
A_{(\dualkegela)}
& := &
\bigoplus_{u \in \dualkegela \cap \dualgittera} A_{u}.%
\end{eqnarray*}
\end{lemma}

\begin{proof}
Let $f_{1}, \ldots, f_{r}$ be homogeneous generators of $A$,
and let $u_{i} \in \dualgittera$ denote the degree of $f_{i}$.
Consider the linear
map $F \colon \Z^{r} \to \dualgittera$ sending the $i$-th canonical
basis vector to $u_{i}$.
Then,
$\gamma := F^{-1}(\dualkegela) \cap \Q^{r}_{\ge 0}$ is a 
pointed, polyhedral cone.
Let $H \subset \gamma$ be the Hilbert Basis of
the semigroup $\gamma \cap \Z^{r}$.
Then, $A_{(\dualkegela)}$ is generated by the elements
$f_{1}^{m_1}\cdot \ldots\cdot f_{r}^{m_r}$,
where $(m_1, \ldots, m_r)  \in H$.%
\end{proof}

\begin{proof}[Proof of Theorem~\ref{data2variety}]
The proof is subdivided into six steps.
Successively weakening the hypotheses,
we prove in the first three steps that 
$\til{X}$ and $X$ are in fact varieties, 
that $\til{X}\to X$ is a proper morphism 
and that 
$\dim(\til{X}) = \dim(Y) + \dim(\torusa)$ holds.
In step four, we show that $\torusa$ acts 
effectively with $\pi \colon \til{X} \to Y$ 
as a good quotient.
Step five is devoted to proving assertion~(iii),
and in step six, we prove normality of $\til{X}$
and $X$.

\begin{step}
\label{regular}
Assume that $\kegela^{\vee} \subset \dualgittera_\Q$ 
is a regular cone, i.e.\ it is mapped to $\Q^r_{\geq 0}$ 
under a suitable isomorphism $\dualgittera \cong \Z^r$.
Moreover, assume $u \mapsto \poldiva(u)$ to be 
linear with integral, basepoint free Cartier divisors
$D_i := \poldiva(e_i)$ where $e_1, \dots, e_r$ are the 
primitive generators of $\kegela^{\vee}$.
In this case, $\til{X} \to Y$ is a rank $r$
vector bundle, and we only have to show that
$\Gamma(Y,\mathcal{A})$ is finitely generated 
and that $\til{X} \to X$ is proper. 
\end{step}

If $r=1$ with an ample Cartier divisor $D_1$,
then we are in the classical setup, cf.~Example~\ref{good}. 
For $D_1$ being just basepoint free, we can reduce to the 
classical case by contracting $Y$ via a
morphism $Y\to \blk{Y}$ with connected fibers such that $D_1$ is
the pull back of an ample Cartier divisor. 
For general $r$, we ``coarsen'' the grading of the
$\mathcal{O}_{Y}$-algebra
$\mathcal{A}$: for $u \in \N^{r}$, denote 
$\vert u \vert := u_{1} + \ldots + u_{r}$ and set
$$
\mathcal{B} 
 \; := \; 
\bigoplus_{k \in \N} \mathcal{B}_{k},
\hspace{1em} \text{where}\hspace{0.7em}
\mathcal{B}_{k} 
\; := \; 
\bigoplus_{\vert u \vert = k} \mathcal{A}_{u}.
$$ 

Consider the corresponding projective space bundle
$Y' := \Proj_Y(\mathcal{B})$ with its projection
$\varphi \colon Y' \to Y$ 
and $\mathcal{L}' := \mathcal{O}_{Y'}(1)$, which means
that $\varphi_*({\mathcal{L}'}^{\otimes k})=\mathcal{B}_k$, 
compare~\cite[p.~162]{Hart}.
Then we obtain an $\mathcal{O}_{Y'}$-algebra and an associated
variety:
$$
\mathcal{A}' \; := \; \bigoplus_{k\in\N}{\mathcal{L}'}^{\otimes k},
\qquad
\til{X}' \; := \; \Spec_{Y'}(\mathcal{A}').$$

Note that $\til{X}'$ is obtained from the rank $r$ vector bundle 
$\til{X}$ over $Y$ by blowing up the zero section 
$s_0 \colon Y \hookrightarrow \til{X}$.
In summary, everything fits nicely into the following commutative 
diagram:
$$
\xymatrix@!0{
& 
{\til{X}'} \ar[rr]  \ar[dl] \ar'[d][dd]
& & 
X' \ar[dd] \ar@{=}[dl]
\\
{\til{X}} \ar[rr]  \ar[dd]
& & 
X  \ar[dd]
\\
&
Y' \ar[dl] \ar'[r][rr] 
& & 
Y'_{0} \ar@{=}[dl]
\\
Y \ar[rr]  \ar@/^1pc/[uu]^-{s_0}
& & 
Y_{0} 
}
$$

The equality $X'=X$ of the spectra of the 
respective rings of global sections follows from
$$ 
\bigoplus_{k \in \N} 
\Gamma(Y', {\mathcal{L}'}^{\otimes k})
\; = \; 
\bigoplus_{k \in \N} 
\Gamma(Y, \varphi_*{\mathcal{L}'}^{\otimes k})
\; = \; 
\bigoplus_{k \in \N} 
\Gamma(Y, \mathcal{B}_k)
\; = \; 
\bigoplus_{u \in \N^r} 
\Gamma(Y, \mathcal{A}_u).
$$

In order to reduce our problem to the case $r=1$,
we have to ascertain that $\mathcal{L}'$ is basepoint free.
Since $\pi(\til{X}'_{f}) = Y'_{f}$ holds
for all homogeneous $f \in \Gamma(Y',\mathcal{A}')$,
it suffices to show that 
any given $x \in \til{X}'\setminus Y'$
admits such an $f$ 
of degree $n \in \N_{>0}$ with $f(x) \ne 0$.
For the latter, consider the canonical projections
$$
\til{X} \; \to \; \til{X}_{i} := \Spec_{Y}(\mathcal{A}_{i}),
\quad \text{where }
\mathcal{A}_{i} \; := \; \bigoplus_{m \in \N} \mathcal{O}(mD_{i}).
$$
Since $\til{X}'\setminus Y'$ equals $\til{X}\setminus Y$,
at least one of these maps
does not send $x \in \til{X} \setminus Y$ to the zero
section of $\til{X}_{i}$. 
By semiampleness of $D_{i}$, there is 
a homogeneous $f \in \Gamma(Y,\mathcal{A}_{i})$ of nontrivial degree 
$me_{i}$ with $f(x) \ne 0$.

\begin{step}\label{finite}
Assume that $\kegela^{\vee} \subset M_\Q$ is a simplicial 
cone, i.e.\ it is generated by linearly independent vectors, 
and that $u \mapsto \poldiva(u)$ is still a linear map.
\end{step}

We will first show that $\til{X}$ is a variety.
For this we only have to verify that $\mathcal{A}$ is locally
of finite type over $\mathcal{O}_Y$.
Choose a sublattice $\unterdualgittera \subset \dualgittera$ 
of finite index such that 
$\unterdualgittera$, $\kegela^{\vee}$ and 
$ u \mapsto \poldiva(u)$ restricted to 
$\kegela^{\vee} \cap \unterdualgittera$
match the assumptions of the previous step.
By linearity, we may extend the assignment
$u \mapsto \poldiva(u)$ from 
$\kegela^{\vee} \cap \dualgittera$ to $\dualgittera$.
This gives further $\mathcal{O}_Y$-algebras 
$$
\mathcal{A}^{\mbox{\rm\tiny grp}}_\unterdualgittera
\;  := \;   
\bigoplus_{u\in \unterdualgittera} \mathcal{A}_{u},
\qquad 
\mathcal{A}^{\mbox{\rm\tiny grp}}_\dualgittera
\;  := \; 
\bigoplus_{u\in \dualgittera} \mathcal{A}_{u}.
$$

Locally, $\mathcal{A}^{\mbox{\rm\tiny grp}}_\unterdualgittera$
looks like $\mathcal{O}_Y\otimes \K[\unterdualgittera]$; hence,
it is locally of finite type.
Choosing representatives $u_1,\dots,u_k$ of 
$\dualgittera/\unterdualgittera$ and
finitely many local generators
$g_{ij}\in\mathcal{A}_{u_i}$,
we see that 
$\mathcal{A}^{\mbox{\rm\tiny grp}}_{\dualgittera}$
is locally finitely generated as an
$\mathcal{A}^{\mbox{\rm\tiny grp}}_\unterdualgittera$-module;
hence, it also is locally of finite type as an
$\mathcal{O}_Y$-algebra.
Finally, we notice that the inclusion
$\mathcal{A}\subset\mathcal{A}^{\mbox{\rm\tiny grp}}_{\dualgittera}$
fits exactly into the situation of 
Lemma~\ref{klaustrick}.
Hence, $\mathcal{A}$ is a locally finitely generated 
$\mathcal{O}_Y$-algebra.

Write for the moment $\mathcal{A}_\dualgittera := \mathcal{A}$,
and, analogously, let $\mathcal{A}_\unterdualgittera$ denote 
the $\mathcal{O}_Y$-algebra associated to
$\poldiva$ restricted to 
$\kegela^{\vee} \cap \unterdualgittera$.
Then the canonical morphism 
$\til{X}_{\dualgittera}\to \til{X}_{\unterdualgittera}$
of the corresponding relative spectra is a finite map 
of varieties; in fact, it is the quotient for the action
of the finite group $\Spec(\K[\dualgittera/\unterdualgittera])$ 
on $\til{X}$ given by the grading.

{From} the preceding step, we know that
$\til{X}_{\unterdualgittera}$ is proper over the affine variety
$X_\unterdualgittera := \Spec(\Gamma(Y,\mathcal{A}_\unterdualgittera))$.
Thus, the affine scheme
$X_\dualgittera := \Spec(\Gamma(\til{X}_\dualgittera,\mathcal{O}))$
gives a commutative diagram
$$
\xymatrix{
{\til{X}}_{\dualgittera} 
\ar[d]_{\text{proper}} 
\ar[rr]^{\text{finite}} 
& & 
{\til{X}}_{\unterdualgittera} \ar[d]^{\text{proper}}   \\
X_{\dualgittera} \ar[rr]_{\text{finite}} 
& & X_\unterdualgittera }
$$
where the lower row is finite because
$\kappa\colon \til{X}_\dualgittera \to X_\unterdualgittera$
is proper, and thus, 
$\Gamma(X_{\dualgittera},\mathcal{O})
= \Gamma(\til{X}_{\dualgittera},\mathcal{O})
= \Gamma(X_{\unterdualgittera},
\kappa_*(\mathcal{O}_{\til{X}_{\dualgittera}}))$
is finite over $\Gamma(X_{\unterdualgittera},\mathcal{O})$.

\begin{step}\label{gluefan}
Let $\poldiva$ be general.
\end{step}

We may subdivide $\kegela^{\vee}$ by a simplicial fan 
$\gkfaecher$ 
such that $\poldiva$ is linear on each of the maximal
cones $\gkfkegela_1,\dots,\gkfkegela_s$ of $\gkfaecher$.
Then, by Step \ref{finite}, 
we know about the corresponding proper maps 
$\til{X}_i\to X_i$. 
The embedding of the cones into the fan yields birational
projections $\til{X}\to\til{X}_i$ and $X\to X_i$,
which in turn lead to closed embeddings
$$
\til{X}\hookrightarrow\til{X}_1\times_Y\ldots\times_Y\til{X}_s,
\qquad
X\hookrightarrow X_1\times_{Y_0}\ldots\times_{Y_0}X_s. $$

\begin{step}\label{gooddim}
The grading of $\mathcal{A}$ defines an effective torus
action $\torusa \times \til{X} \to \til{X}$
having the canonical map $\til{X} \to Y$ as a
good quotient. 
\end{step}

For any affine $V \subset Y$, the grading of 
$\Gamma(V,\mathcal{A})$
defines a $\torusa$-action on $\Spec(\Gamma(V,\mathcal{A}))$.
This is compatible with glueing, and thus we obtain a
$\torusa$-action on $\til{X}$.
By construction, $\til{X} \to Y$ is affine, and
$\mathcal{O}_{Y}=\mathcal{A}_0$
is the sheaf of invariants.
Hence, $\til{X} \to Y$ is a good quotient for the $\torusa$-action.
The fact that $\torusa$ acts effectively is seen as follows.
Since the algebra $\mathcal{A}$ admits locally nontrivial 
sections in any degree $u \in \kegela^{\vee} \cap \dualgittera$,
the weight monoid of the $\torusa$-action generates
$\dualgittera$ as a group.
Consequently, the $\torusa$-action has free orbits, and hence is
effective.%

\begin{step}
\label{part3}
Let $f \in A_u$. Here we will prove the third part 
of the theorem.
Note that by Condition \ref{cpldef}~(iii) of the map
$u \mapsto \poldiva(u)$,
this will imply birationality of $\til{X} \to X$.
\end{step}

In the situation of Step~\ref{regular},
the equality $\quotabb(\til{X}_{f}) = Y_{f}$
is obvious (and was already used there).
Suppose we are in the setting of Step~\ref{finite}.
There we introduced a finite map
$\til{X}_\dualgittera \to 
\til{X}_\unterdualgittera$ 
where $\til{X}_\dualgittera = \til{X}$.
This map fits into the commutative diagram 
$$ 
\xymatrix{
{\til{X}_\dualgittera}
\ar[rr] \ar[dr]_{\pi_\dualgittera}
& & 
{\til{X}_\unterdualgittera}
\ar[dl]^{\pi_\unterdualgittera}
\\
& Y & 
}
$$
where we denote by 
$\pi_\dualgittera$ and 
$\pi_\unterdualgittera$ the respective
canonical projections.
Then $f^k \in \Gamma(\til{X}_\unterdualgittera,\mathcal{O})$ 
holds for some positive power of 
$f \in \Gamma(\til{X}_\dualgittera,\mathcal{O})$,
and, clearly,  
$\pi_\dualgittera((\til{X}_\dualgittera)_f)$
equals
$\pi_\unterdualgittera((\til{X}_\unterdualgittera)_{f^k})$.
This reduces the problem to the setting
of Step~\ref{regular}.

Thus, we are left with considering the situation of 
Step~\ref{gluefan}.
There we used a simplicial fan subdivision of $\kegela^{\vee}$
with maximal cones $\gkfkegela_i$.
This defines birational morphisms 
$\varphi_{i} \colon \til{X} \to \til{X}_{i}$
and commutative diagrams

$$
\xymatrix{
{\til{X}} \ar[rr]^{\varphi_{i}}\ar[dr]_{\quotabb}
& & 
{\til{X}_{i}} \ar[dl]^{\quotabb_{i}} \\
& Y & }
$$

Choose $i$ such that $u = \deg(f)$
lies in $\gkfkegela_{i}$, and
write $f = \varphi_{i}^{*}(f_{i})$.
Then the diagram directly gives 
$\quotabb(\til{X}_f) \subset Y_f$.
For the converse inclusion,
note first that $\varphi_{i}$ 
induces dominant morphisms 
of the fibers $\til{X}_{y} \to \til{X}_{i,y}$.
Now, let $y \in Y$ such that $f$
vanishes along $\quotabb^{-1}(y)$.
Then, by fiber-wise dominance of $\varphi_{i}$,
the function $f_{i}$ vanishes along $\quotabb_{i}^{-1}(y)$.
Hence, the previous cases result in $y \not\in Y_{f}$.

If $Y_{f}$ is affine,
then so is $\quotabb^{-1}(Y_{f})$ and hence, by  
$\til{X}_{f} \subset \quotabb^{-1}(Y_{f})$,
also $\til{X}_{f} = \quotabb^{-1}(Y_{f})_{f}$.
In particular, we have an isomorphism
 $\til{X}_{f} \to X_{f}$ in this case.
The last statement of (iii) can be proven as follows:
for any $v \in \kegela^{\vee} \cap \dualgittera$, we have
$$
\begin{array}{rcl}
\Gamma(Y_f,\mathcal{A}_v)&=&
\{g\in \K(Y)\kst \ddiv(g) + \poldiva(v)
\geq 0\text{ on } Y_f\}\\
&=&
\{g\in \K(Y)\kst \exists k\geq 0:
\ddiv(g) + \poldiva(v) + k\cdot(\ddiv(f) + \poldiva(u))\geq 0\}\\
&=&
\{g\in \K(Y)\kst \exists k\geq 0:
gf^k\in A_{v+ku}\} \\
& = & (A_f)_v.
\end{array}
$$

\begin{step}\label{normal}
The varieties $\til{X}$ and $X$ are normal.
\end{step}

It suffices to show that $\til{X}$ is normal. 
This is a local problem; hence, we may assume 
in this step that $Y$ is affine 
and, moreover, 
that for all $u \in \kegela^{\vee} \cap \dualgittera$,
the homogeneous pieces 
$A_u := \Gamma(Y, \mathcal{A}_u)$ of $A$ 
are non-trivial.
We may use any 
$g_{u}\in A_{u}$ to obtain an embedding
$$
\imath_{u} \colon A_{u} \; \hookrightarrow \; Q(A)^{\torusa} = Q(A_{0}) =
\K(Y), 
\qquad 
f \; \mapsto \; f/g_{u},
$$
where $Q(A)$ stands for the fraction field of $A$.
Note that the equality $Q(A)^{\torusa} = Q(A_0)$ 
holds because the quotient space
$Y$ is of dimension $\dim(\til{X}) - \dim(\torusa)$. 
The image of the above embedding can be described as follows:
\begin{eqnarray*}
\imath_{u}(A_{u})
& = &
\{ h \in Q(A_{0}) \kst
      \ddiv(h) \ge - \poldiva(u) - \ddiv(g_{u}) \}.
\end{eqnarray*}

Now, we consider the integral closure $\blkk{A}$ of $A$ in
$Q(A)$.
Recall that $\blkk{A}$ is also $\dualgittera$-graded, 
see e.g.~\cite[Prop.~V.8.21]{Bou}. 
Thus, in order to show $A = \blkk{A}$, we only 
have to verify that a homogeneous $f \in Q(A)$,
say of degree $u \in M$, 
belongs to $A$ provided it satisfies a 
homogeneous equation
of integral dependence with certain $h_{i} \in A$:
\begin{eqnarray*}
f^{n} 
& = & 
h_{1}f^{n-1} + \ldots + h_{n-1}f + h_{n}. 
\end{eqnarray*}

This equation implies, in particular, 
that the degree $u\in\dualgittera$ belongs to
the weight cone $\kegela^{\vee}$.
Hence, we may choose an element $g_{u}\in A_{u}$.
Suppose, for the moment, that $\poldiva(u)$ is an 
integral Cartier divisor with its sheaf being
locally generated, 
without loss of generality, by $g_u$.
Then the above equation expressing the integral dependence
of $f$ takes place over
$$
B 
\; := \: 
\bigoplus_{n \in \N} A_{nu}
\; = \: 
\bigoplus_{n \in \N} \Gamma(Y, \mathcal{O}(n\poldiva(u)))
\; = \:
A_0[g_u],
$$
and we are done because of 
$f\in g_uQ(A_0)\subset Q(B)$ and
the integral closedness of $B$.
In the general case, we choose an $m > 0$ such that $m\poldiva(u)$ 
is an integral Cartier divisor.
The previous argument yields
$f^{m} \in A_{mu}$. Hence,
by the description of $A_{mu}$ in terms of 
the injection $\imath_{mu} \colon A_{mu} \to Q(A_{0})$, 
this means that
\begin{eqnarray*} 
\ddiv(f^{m}/g_{u}^{m}) 
& \ge & 
- m \,\poldiva(u) - m \ddiv(g_{u}). 
\end{eqnarray*}
Dividing this inequality by $m$ shows 
$f/g_{u}\in \imath_{u}(A_{u})$.
This in turn means $f \in A_{u}$, which concludes the proof.
\end{proof}
\setcounter{step}{0}

\section{Ingredients from GIT}

In this section, we recall crucial ingredients 
from Geometric Invariant Theory 
for the proof of Theorem~\ref{variety2data} 
and also for the applications  
presented later.
The central statement is a 
description 
of the GIT-equivalence classes arising 
from  linearizations
of the trivial bundle over an affine 
variety with torus action in terms of a quasifan.
For torus actions on $\K^n$, this 
is well known; the describing quasifan
then is even a fan, and is called
a Gelfand-Kapranov-Zelevinsky decomposition,
compare~\cite{OdPa}.
For details on the general case as presented 
here, we refer to~\cite{BeHa}.

Let us remark that there are analogous,
and even further going results in the 
projective case.
Brion and Procesi~\cite{BriPro}
observed that, for a torus action on a projective
variety, the collection of all GIT-quotients
arising from the different linearizations of a 
given ample bundle
comes along with a piecewise linear structure.
This has been generalized in~\cite{Thadd},
~\cite{DoHu} and~\cite{Re}
to arbitrary reductive groups and the collection of all 
GIT-quotients arising from linearized ample bundles;
see also~\cite{Chow} for some work in the 
toric setup.

Let us fix the setup.
By $\dualgittera$, we denote a lattice, 
and $A$ is an integral, affine,
$\dualgittera$-graded $\K$-algebra:
\begin{eqnarray*}
A
& = & 
\bigoplus_{u \in \dualgittera} 
A_u.
\end{eqnarray*}
Let $X := \Spec(A)$ denote the affine 
variety associated to $A$.
Then the $\dualgittera$-grading of 
$A$ defines an action of the algebraic
torus $T := \Spec(\K[\dualgittera])$ 
on $X$.

For convenience, we briefly recall the 
basic concepts from~\cite{GIT} in a 
down-to-earth manner.
A {\em $T$-linearization\/} of a line bundle
$L \to X$ is a fiberwise linear $T$-action 
on the total space $L$ such that the bundle 
projection $L \to X$ becomes $T$-equivariant.
Any $T$-linearization of the trivial 
bundle over $X$ is of the form 
\begin{eqnarray}
\label{lintriv}
t \mal (x,z)
& = &
(t \mal x, \chi^u(t)z), 
\end{eqnarray}
where $\chi^u \colon T \to \K^*$ denotes 
the character corresponding to $u \in \dualgittera$.
Note that the $n$-fold tensor product of the above 
linearization corresponds to the character
$\chi^{nu}$.
Any $T$-linearization of a line bundle defines
a representation of $T$ on the space of its 
sections via
\begin{eqnarray*}
(t \mal s)(x)
& := & 
t \mal (s (t^{-1} \mal x)).
\end{eqnarray*}

The set of {\em semistable points\/} 
associated to a $T$-linearized line bundle $L \to X$ 
is defined as the union of all 
affine sets of the form $X_f$, where $f$ is a 
$T$-invariant section of some positive 
tensor power $L^{\otimes n}$. 
The invariant sections for the 
linearization~\ref{lintriv} are precisely
the elements $f \in A_{nu}$, where $n \in \Z_{> 0}$,
and the corresponding set of semistable points
is
\begin{eqnarray*}
X^{ss}(u)
& := & 
\bigcup_{f \in A_{nu},  \ n \in \Z_{> 0}} X_{f}.
\end{eqnarray*}

Two linearized bundles are called 
{\em GIT-equivalent\/} if they define the 
same sets of semistable points.
The description of the GIT-equivalence classes 
arising from the linearizations of the trivial 
bundle presented in~\cite{BeHa} works in
terms of orbit cones.
Let us recall the definition of these
and other orbit data.

\begin{definition}
\label{orbitdatadef}
Consider a point $x \in X$.
\begin{enumerate}
\item
The {\em orbit monoid\/} associated to
$x \in X$ is the submonoid $S(x) \subset M$ 
consisting of all $u \in \dualgittera$ 
that admit an $f \in A_u$ with $f(x) \ne 0$
\item 
The {\em orbit cone\/} associated to $x \in X$ is 
the convex cone $\omega(x) \subset \dualgittera_\Q$ 
generated by the orbit monoid.
\item 
The {\em orbit lattice\/} associated to
$x \in X$ is the sublattice $M(x) \subset M$ 
generated by the orbit monoid. 
\end{enumerate}
\end{definition}

The orbit cones are polyhedral, and each of them
is contained in the weight cone 
$\omega \subset \dualgittera_\Q$,
which in turn is generated by the $u \in \dualgittera$ 
with $A_u \ne \{0\}$.
The geometric meaning of the above orbit data
is made clear by the following:

\begin{proposition}
\label{orbitcones}
Consider a point $x \in X$.
\begin{enumerate}
\item 
The orbit lattice $M(x)$ consits of all
$u \in M$ which admit a $u$-homogeneous function 
$f \in \K(X)$ that is defined and invertible 
near $x$.
\item 
The isotropy group $T_x \subset T$ of 
the point $x \in X$ is the diagonalizable 
group given by $T_x = \Spec(\K[M/M(x)])$.
\item 
The orbit closure $\blk{T \mal x}$ is 
isomorphic to $\Spec(\K[S(x)])$; it comes 
along with an equivariant open embedding
of the torus $T/T_x = \Spec(\K[M(x)])$.
\item 
The normalization of the orbit closure
$\blk{T \mal x}$ is the toric variety corresponding 
to the cone $\omega(x)^{\vee}$ in $\Hom(M(x),\Z)$.
\end{enumerate}
\end{proposition}

In terms of orbit cones, there is a simple description
of the sets $X^{ss}(u)$ of semistable points.
Namely, we have
\begin{eqnarray*}
X^{ss}(u)
& = & 
\{x \in X; \; u \in \omega(x)\}.
\end{eqnarray*}

\begin{definition}
The {\em GIT-cone\/} associated to an element 
$u \in \omega \cap M$ 
is the intersection of all orbit cones 
containing $u$:
\begin{eqnarray*}
\gkfkegela(u)
& := & 
\bigcap_{x \in X; \; u \in \omega(x)}
\omega(x).
\end{eqnarray*}
\end{definition}

The GIT-cones turn out to be polyhedral cones
as well. Their importance is that they
correspond to the GIT-equivalence classes. 
The main results of~\cite[Section~2]{BeHa}
can be summarized as follows:

\begin{theorem}
\label{gitfan}
Let $A$ be an integral affine algebra
graded by a lattice $\dualgittera$.
Then, for the action of $\torusa := \Spec(\K[\dualgittera])$
on $X := \Spec(A)$, the following statements hold. 
\begin{enumerate}
\item
The GIT-cones $\gkfkegela(u)$, where $u \in \dualgittera$,
form a quasifan $\gkfaecher$ in $\dualgittera_\Q$.
\item The support of the quasifan  $\gkfaecher$ is the
weight cone $\omega \subset \dualgittera_\Q$.
\item
For any two elements $u_1, u_2 \in \omega \cap \dualgittera$,
we have
\begin{eqnarray*}
X^{ss}(u_1) \subset X^{ss}(u_2)
& \iff &
\gkfkegela(u_1) \supset \gkfkegela(u_2).
\end{eqnarray*}
\end{enumerate}
\end{theorem}

The set of semistable points
of a $T$-linearized line bundle over $X$
is an open $T$-invariant subset of $X$, 
and it admits a good quotient by the action
of $T$. 
For the linearization~\ref{lintriv},
the quotient space
$Y_u := X^{ss}(u) \quot T$ is given by 
$$ 
Y_u
\; = \; 
\Proj\left(A_{(u)}\right),
\qquad
\text{where} \quad
A_{(u)} 
\; = \; \bigoplus_{n \in \Z_{\ge 0}} A_{nu}
\; \subset \; 
A.  
$$
Note that each quotient space $Y_u$ is projective 
over the affine variety $Y_0 = \Spec(A_0)$. 
For our purposes, the following observation concerning
the dimension of quotient spaces will be
needed.

\begin{lemma} 
Suppose that the action of $T$ on $X$ is effective.
Then, for every $u \in \relint{\dualkegela} \cap \dualgittera$, 
we have
\begin{eqnarray*}
\dim(X^{ss}(u) \quot T)
& = & 
\dim(X) - \dim(\torusa).
\end{eqnarray*}
\end{lemma}

\begin{proof}
Note that the orbit cone of a generic
orbit $T \mal x_0 \subset X$
equals the weight cone $\omega$.
Thus, for $u \in \relint{\omega}$,
the orbit $T \mal x_0$ is a 
closed subset of $X^{ss}(u)$.
Thus, $T \mal x_0$ is a generic
fiber of 
$X^{ss}(u) \to X^{ss}(u) \quot T$. 
Since the $T$-action is effective,
we have $\dim(T \mal x_0) = \dim(T)$.
The assertion follows.
\end{proof}

\section{Proof of Theorem~\ref{variety2data}}
\label{limit}

In this section, we prove Theorem~\ref{variety2data}.
The setup is the following:
$\dualgittera$ is a lattice, and $A$ is a 
$\dualgittera$-graded affine $\K$-algebra.
We consider $X := \Spec(A)$ and the action 
of $\torusa := \Spec(\K[\dualgittera])$
on $X$ defined by the grading;
we assume that $X$ is normal and that the
$\torusa$-action is effective.
By $\gewkegel \subset \dualgittera_{\Q}$, 
we denote the weight cone of the $\torusa$-action, 
and $\gkfaecher$ is the quasifan in $\dualgittera_\Q$ 
consisting of the GIT-cones as 
discussed in Theorem~\ref{gitfan}.

Then, for every $\gkfkegela \in \gkfaecher$,  the map 
$u \mapsto X^{ss}(u)$ 
is constant on the relative interior $\relint{\gkfkegela}$.
We denote by $W_{\gkfkegela} \subset X$ the set of 
semistable points
defined by any of those $u \in \relint{\gkfkegela}$
and, moreover, by
$q_{\gkfkegela} \colon W_{\gkfkegela} \to Y_{\gkfkegela}$,
the corresponding good $\torusa$-quotient.
In particular, we have $W_{0} = X= \Spec(A)$
and $Y_0=\Spec(A_0)$.
The spaces $Y_{\gkfkegela}$ are normal; 
the morphisms $q_{\gkfkegela}$ are affine, 
and each of their fibers contains exactly
one closed $\torusa$-orbit, and hence it is connected.

The quotient maps
$q_{\gkfkegela} \colon W_{\gkfkegela} \to Y_{\gkfkegela}$,
where $\gkfkegela \in \gkfaecher$,
can be put together to an inverse system with 
$q_{0} \colon W_0 \to Y_0$ sitting at the end.
Let us consider its inverse limit.
If $\gkfkegelb \preceq \gkfkegela$,
then we have an open embedding
$W_\gkfkegela \subset W_\gkfkegelb$ 
inside $X$. Set
$$
\begin{array}{ccccc}
W 
& := & 
\kprojlim 
W_\gkfkegela
& = & 
\bigcap_{\gkfkegela\in\gkfaecher}W_{\gkfkegela}.
\end{array}
$$
The inverse limit $Y'$ of the induced maps
$p_{\gkfkegela\gkfkegelb} \colon Y_\gkfkegela\to Y_\gkfkegelb$
between the quotient spaces
is a nested fiber product 
over $Y_0$.
The inverse limit of all quotient maps 
is the canonical map
$q' \colon W \to Y'$.

In general, $Y'$ might be reducible,
but there is a canonical irreducible 
component:
the closure of the image $q'(W)$. 
We obtain an irreducible, normal variety 
$Y$ by
taking the normalization of this canonical
component:
\begin{eqnarray*}
Y
& := & \text{normalization}
\left(
\blk{q'(W)}
\right).
\end{eqnarray*}
By the universal property of the normalization,
we have an induced morphism 
$q \colon W \to Y$. In summary, we obtain 
for each pair $\gkfkegelb \preceq \gkfkegela$
in $\gkfaecher$ a commutative diagram
$$
\xymatrix{
W \ar[rr]^-{j_\gkfkegela} \ar[d]_{q} & & 
W_{\gkfkegela} \ar[d]^{q_{\gkfkegela}} \ar[rr]^-{j_{\gkfkegela\gkfkegelb}} 
& & W_{\gkfkegelb} \ar[d]^-{q_{\gkfkegelb}} 
\ar[r]^-{j_{\gkfkegelb 0}} & X \ar[dd]^-{q_0}\\
Y \ar[rr]^-{p_{\gkfkegela}} \ar[drrrrr]_-{p_0} & & Y_{\gkfkegela} 
\ar[drrr]^-{p_{\gkfkegela 0}} \ar[rr]^-{p_{\gkfkegela\gkfkegelb}}
& & Y_{\gkfkegelb}\ar[dr]^-{p_{\gkfkegelb 0}} \\
&&&&& Y_{0} & 
}
$$

\begin{lemma}\label{projproj}
The morphisms $p_{\gkfkegela} \colon Y \to Y_{\gkfkegela}$ 
and 
$p_{\gkfkegela\gkfkegelb} \colon Y_{\gkfkegela} \to Y_{\gkfkegelb}$
are projective 
surjections with connected fibers.
Moreover, if $\dim(Y_{\gkfkegela})=\dim(X)-\dim(\torusa)$,
for example if $\gkfkegela$ intersects $\relint{\gewkegel}$,
then the morphism 
$p_{\gkfkegela} \colon Y \to Y_{\gkfkegela}$ is birational.
\end{lemma}

\begin{proof}
First, recall 
that each quotient space $Y_{\gkfkegela}$ 
is projective over $Y_{0}$. 
It follows that $Y$ is projective over $Y_{0}$, and thus, that 
each of the maps $p_{\gkfkegela}$ is projective, too. 
Since every $Y_{\gkfkegela}$ is dominated by $W$, all morphisms
$p_{\gkfkegela} \colon Y \to Y_{\gkfkegela}$ are dominant.
Together with properness, this implies surjectivity of
each  $p_{\gkfkegela}$.
The same reasoning leads to these properties for the
$p_{\gkfkegela\gkfkegelb}$.

Let us show the connectedness of the fibers.
In general, the fibers of a projective morphism 
between normal varieties are connected 
if and only if its generic fiber is is connected;
use, e.g., Stein factorization.
Thus, it suffices to check that  
the generic fiber of 
$p_{\gkfkegela} \colon Y \to Y_{\gkfkegela}$
is irreducible.

%
The image $q(W) \subset Y$ is constructible. Hence,
we can choose a closed proper subset $C \subset Y$
such that $q(W)$ and $C$ cover $Y$.
Let $y \in Y_{\gkfkegela}$ be a generic point.
Then, the  fiber $p_{\gkfkegela}^{-1}(y) \subset Y$ splits into
$$ 
\begin{array}{ccccc}
p_{\gkfkegela}^{-1}(y)
& = &
\blk{q (q^{-1}(p_{\gkfkegela}^{-1}(y)))} 
\cup (p_{\gkfkegela}^{-1}(y) \cap C)
& = & 
\blk{q (q_{\gkfkegela}^{-1}(y)\cap W)} 
\cup (p_{\gkfkegela}^{-1}(y) \cap C).
\end{array}
$$
Since we know that the generic fiber of
$q_{\gkfkegela} \colon W_{\gkfkegela} \to Y_{\gkfkegela}$
is irreducible~\cite[Prop.~4]{Ar}, 
this also holds for the first part of the previous
expression.
Thus, it suffices to show that
this part actually fills
the whole fiber $p_{\gkfkegela}^{-1}(y)$.

Assume to the contrary that there is an
irreducible component 
$C_{0} \subset C$ 
dominating $Y_{\gkfkegela}$ and containing
some irreducible component $F_{0} \subset p_{\gkfkegela}^{-1}(y)$.
By general properties of morphisms,
see for example~\cite[II.3.22]{Hart}
or~\cite[I.4.1 and I.4.3]{Humph},
this would imply
$$
\begin{array}{ccccc}
\dim(F_{0})
&  \geq &  
\dim(Y) - \dim(Y_{\gkfkegela})
& > & 
\dim(C_{0}) - \dim(Y_{\gkfkegela}).
\end{array}
$$
On the other hand, consider the restriction 
$\pi_{\gkfkegela} \colon C_{0} \to Y_{\gkfkegela}$ of $p_{\gkfkegela}$.
Since $y\in Y_{\gkfkegela}$ is generic, the dimension
of $\pi_{\gkfkegela}^{-1}(y)$ equals $\dim(C_{0}) - \dim(Y_{\gkfkegela})$.
This contradicts the previous estimation.

Finally, we need to prove
the claim about the maps
$p_\gkfkegela$ being birational.
If we are given two cones $\gkfkegelb\preceq \gkfkegela$ both
satisfying the assumption
$\dim(Y_{\gkfkegelb})=\dim(Y_{\gkfkegela})=\dim(X)-\dim(\torusa)$,
then the connecting map 
$p_{\lambda \gamma} \colon Y_\lambda \to Y_\gamma$
is birational, because it induces the identity map
on
$$ 
\K(Y_\gamma)
\; = \; 
\K(Y_\lambda)
\; = \; 
\K(X)^T.
$$

But $Y$ can also be obtained from the
complete inverse system provided by all cones
$\gkfkegela\in\gkfaecher$
which intersect $\relint{\gewkegel}$.
Thus, $Y$ can be built from a system of birational maps,
and the common open subset (where all the
$p_{\gkfkegela\gkfkegelb}$ are isomorphisms) survives in $Y$.
\end{proof}

We will now investigate certain coherent sheaves on the
quotient spaces $Y_{\gkfkegela}$.
As mentioned earlier, we have $Y_{\gkfkegela} = \Proj(A_{(u)})$ 
with the ring
\begin{eqnarray*}
A_{(u)}
& = &
\bigoplus_{n \in \N} A_{nu},
\end{eqnarray*}
where $u \in \relint{\gkfkegela} \cap M$ may 
be any element.
This allows us to associate to $u$ a sheaf on 
$Y_{\gkfkegela}$, namely
$$
\mathcal{A}_{\gkfkegela,u} 
\; := \; 
\mathcal{O}_{Y_\gkfkegela}(1)
\; = \;
(q_{\gkfkegela})_{*}(\mathcal{O}_{W_{\gkfkegela}})_{u},
$$
where in the last expression, the subscript ``$u$''
indicates that we mean the sheaf of 
semiinvariants with respect to the character
$\chi^u \colon \torusa \to \dualgittera$.

\begin{remark}   
In the terminology of~\cite{Hart}, our 
$\mathcal{A}_{\gkfkegela,u}$ is nothing but
the sheaf on the $\Proj$ associated to the graded 
module $A_{(u)}(1)$.
\end{remark}


We call an element $u \in \dualgittera$ {\em saturated} 
if the ring $A_{(u)}$ is generated in degree one.
It is well known~\cite[Prop.~II.1.3]{Bou} 
that every $u \in \dualgittera$ admits a positive 
multiple $nu \in \dualgittera$ 
such that all positive multiples of $nu$ are 
saturated.
Moreover, as before,
denote by $Q(A)$ the field of fractions of $A$.

\begin{lemma}
\label{amplesheaves}
Let $\gkfkegela \in \gkfaecher$ and $u \in \relint{\gkfkegela}$.
For $f \in A_{nu}$, let
$Y_{\gkfkegela,f} := q_{\gkfkegela}(X_{f})$
be the corresponding affine chart of $Y_{\gkfkegela}$.
\begin{enumerate}
\item On $Y_{\gkfkegela,f}$, the sheaf $\mathcal{A}_{\gkfkegela,u}$
  is the coherent $\mathcal{O}_{Y_{\gkfkegela}}$-module associated to the
  $(A_{f})_{0}$-module $(A_{f})_{u}$.
\item
  If $u$ is saturated, then $\mathcal{A}_{\gkfkegela,u}$ 
  is an ample invertible sheaf on $Y_{\gkfkegela}$, and
  on the charts $Y_{\gkfkegela,f}$, where $f \in A_{u}$, 
  we have
  $$
  \mathcal{A}_{\gkfkegela,u}
  \; = \; 
  f\cdot(A_f)_0
  \; = \; 
  f\cdot \mathcal{O}_{Y_\gkfkegela}.
  $$
\item If $g \in Q(A)$ and $n \in \Z_{>0}$, then
  $g^n \in \mathcal{A}_{\gkfkegela,nu}$ implies
  $g \in \mathcal{A}_{\gkfkegela,u}$.
\item 
  The global sections of $\mathcal{A}_{\gkfkegela,u}$ are
  $\Gamma(Y_\gkfkegela, \mathcal{A}_{\gkfkegela,u})=A_u$.
\end{enumerate}
\end{lemma}

\begin{proof} 
The first two assertions are obvious.
To prove the third one,
let $g\in Q(A)$ such that $g^n\in (A_f)_{nu}$. 
We may assume that $f$ appears in the
denominator with a power divisible by $n$. Then,
there is some $k\geq 0$ such that $(gf^k)^n\in A$.
The normality of $A$ implies $gf^k\in A$; thus,
$g\in (A_f)_{u}$.

We turn to the last statement. 
Let $u \in \relint{\gkfkegela} \cap M$. 
Since $\Gamma(Y_{\gkfkegela},\mathcal{A}_{\gkfkegela,u})$ equals 
$\Gamma(W_{\gkfkegela},({\mathcal{O}}_{W_{\gkfkegela}})_{u})$,
we need to prove that any $u$-homogeneous function
$g$ on $W_{\gkfkegela}$ extends to $X$.
By normality of $X$, it suffices to show that $g$ has non-negative
order along any prime divisor contained in $X \setminus W_{\gkfkegela}$.
For the latter, we may also take any positive power $g^{n}$.
Thus, we may assume that $u$ is saturated.

Consider a prime divisor $D \subset X \setminus W_{\gkfkegela}$.
Choose $f \in A_{u}$ such that the order 
$\nu_{D}(f)$ of $f$ along $D$ is minimal.
Regarding $g$ as an element of $(A_{f})_{u}$,
we find a $k \geq 0$ and an $h \in A_{(k+1)u}$ such that
$g = h/f^{k}$.
Since the elements of $A_{(k+1)u}$ are polynomials in elements of
$A_{u}$, the minimality of $\nu_{D}(f)$ provides
$\nu_{D}(h)\geq (k+1)\nu_{D}(f)$; hence
$\nu_{D}(g)\geq \nu_{D}(f)\geq 0$.
\end{proof}

The sheaves $\mathcal{A}_{\gkfkegela,u}$ live on different 
spaces.
By pulling them back, we obtain for every
$u \in \gewkegel \cap \dualgittera$ a well defined
coherent sheaf on $Y$:
$$
\mathcal{A}_u
\; := \; p_\gkfkegela^* (\mathcal{A}_{\gkfkegela,u}),
\text{ where }
\gkfkegela \in \gkfaecher \text{ is the cone with } 
u\in\relint{\gkfkegela}.
$$

\begin{lemma}\label{bpfree}
Let $u, u' \in \gewkegel \cap \dualgittera$.
\begin{enumerate}
\item
We have $\K(Y)=Q(A)_0$, and the natural transformation
$p_{\gkfkegela}^{*}{q_{\gkfkegela}}_{*} \to q_* j_\gkfkegela^*$
sends $\mathcal{A}_u$ into $Q(A)_u$.
\item
Let $u$ be saturated. Then $\mathcal{A}_u$ is a globally generated, 
invertible sheaf. On the (not necessarily affine) sets
$Y_f := p_\gkfkegela^{-1}(Y_{\gkfkegela,f})$
with $f\in A_{u}$, we have
$$\mathcal{A}_u \; = \; f\cdot\mathcal{O}_Y \; 
\subset \; f\cdot \K(Y) \; = \; Q(A)_u.$$
Moreover, for the global sections of $\mathcal{A}_{u}$, we obtain
$\Gamma(Y,\mathcal{A}_{u}) = A_u$.
\item
If $u,u'$ and $u+u'$ are saturated, then we have
$\mathcal{A}_{u}\mathcal{A}_{u'}\subset\mathcal{A}_{u+u'}$. 
If, moreover, $u$ and $u'$ lie in a common cone of 
$\gkfaecher$, then equality holds.
\end{enumerate}
\end{lemma}

\begin{proof}
Using an arbitrary cone 
$\gkfkegela \in \gkfaecher$ which intersects 
$\relint{\gewkegel}$
and some homogeneous $f\in A$ with 
$\deg(f) \in \relint{\gkfkegela}$,
we obtain 
$$\K(Y)\; =\; \K(Y_\gkfkegela)\; =\; Q((A_f)_0)\; \subset \; Q(A)_0.$$
Conversely, starting with an element 
$a/b\in Q(A)_0$, then
$\deg a=\deg b$ is sitting in the interior of
some cone $\gkfkegelb \in \gkfaecher$; thus, we have
$$a/b\; \in \; Q((A_b)_0)\; =\; \K(Y_\gkfkegelb)\; \subset\; \K(Y).$$

The second assertion is a direct consequence
of Lemma \ref{amplesheaves}. For the last part, we use that
the adjunction map 
$\mathcal{A}_{\gkfkegela,u} \to 
p_{\gkfkegela*}p_\gkfkegela^*\mathcal{A}_{\gkfkegela,u}$
locally looks like
$\mathcal{O}_{Y_\gkfkegela}
\to 
p_{\gkfkegela*}p_\gkfkegela^*\mathcal{O}_{Y_\gkfkegela}$.
Hence, because of Lemma \ref{projproj}, it is an isomorphism.

Eventually, to prove the third assertion,
we have to deal with the product
$\mathcal{A}_{u}\mathcal{A}_{u'}$.
Due to saturatedness, we obtain 
$$
\mathcal{A}_{u}\mathcal{A}_{u'}
\; = \; 
(A_u\cdot\mathcal{O}_Y)  (A_{u'}\cdot\mathcal{O}_Y)
\; = \; 
A_u A_{u'} \cdot\mathcal{O}_Y 
\; \subset\; 
A_{u+u'} \cdot\mathcal{O}_Y
\; = \; 
\mathcal{A}_{u+u'}.
$$

Now let $u \in \relint{\gkfkegelb}$, $u' \in \relint{\gkfkegelb}'$
and, moreover, 
$\gkfkegelb,\gkfkegelb' \subset \gkfkegela \in \gkfaecher$
(with $\gkfkegela$ being minimal).
Then we have 
$W_{\gkfkegela} \subset W_{\gkfkegelb} \cap W_{\gkfkegelb'}$.
Conversely, again by saturatedness, we obtain
$$
W_{\gkfkegelb} \cap W_{\gkfkegelb'} 
  =
\biggl(\bigcup_{f \in A_u} X_f \biggr) 
 \cap 
\biggl(\bigcup_{f' \in A_{u'}} X_{f'} \biggr) 
 = 
\bigcup_{\genfrac{}{}{0pt}{}
         {\mbox{$\scriptstyle f\in A_u$}}
         {\mbox{$\scriptstyle f'\in A_{u'}$}}}
      X_{ff'}  
\subset 
\bigcup_{g\in A_{u+u'}} X_g 
 =  W_\gkfkegela.
$$

In particular, the sets $X_{ff'}$ cover $W_\gkfkegela$,
hence, the $Y_{\gkfkegela,ff'}$ cover $Y_\gkfkegela$,
hence, so do the $Y_{ff'}$ with $Y$.
On the other hand, the inclusion
$X_{ff'}\subset X_{f}$ induces a morphism
$Y_{\gkfkegela,ff'}\to Y_{\gkfkegelb,f}$
which also applies to $f'$.
Thus, $Y_{ff'}\subset Y_f\cap Y_{f'}$,
and it follows that 
$\mathcal{A}_u=f\cdot \mathcal{O}_Y$,
$\mathcal{A}_{u'}=f'\cdot \mathcal{O}_Y$ and
$\mathcal{A}_{u+u'}=ff'\cdot \mathcal{O}_Y$
on $Y_{ff'}$.
\end{proof}


\begin{proof}[Proof of Theorem \ref{variety2data}]
Let $Y$ be the semiprojective variety 
defined at the beginning of this section. 
We will construct the desired pp-divisor 
on $Y$ as a convex piecewise linear map  
$\gewkegel \to \CDiv_{\Q}(Y)$
in the sense of Definition~\ref{cpldef}.
The previous two Lemmas will be used implicitly.

Our construction requires a (non-canonical) 
choice of a homomorphism
$s \colon \dualgittera \to Q(A)^{*}$
such that for every $u \in \dualgittera$ the function
$s(u)$ is homogeneous of degree $u$.
Since $\torusa$ acts effectively on $X$, 
such ``sections'' $s \colon \dualgittera \to Q(A)^{*}$
always exist.

Now, if $u \in \gewkegel \cap \dualgittera$ is any 
saturated element, then there is a unique Cartier divisor
$\poldiva(u) \in \CDiv(Y)$ such that
$$
\mathcal{O}_Y\big(\poldiva(u)\big) 
\; = \; 
\frac{1}{s(u)} \cdot \mathcal{A}_u
\; \subset \; \K(Y).
$$
The local equation for 
$\poldiva(u)$ on $Y_f$ with $f\in A_u$
is $s(u)/f$.
For general $u \in \gewkegel$, we can choose a saturated
multiple $nu$ and define
$$
\poldiva(u) \; := \; \frac{1}{n} \cdot \poldiva(nu) \; \in \; \CDiv_{\Q}(Y).
$$

Obviously, this definition does not depend on $n$,
and one directly
checks the properties of Definition~\ref{cpldef}
for the map $u \mapsto \poldiva(u)$.
Moreover, we can recover
the $\dualgittera$-graded ring $A$ via
$$
A_{u} \; = \; s(u) \cdot \Gamma(Y,\mathcal{O}_Y(\poldiva(u))) .
$$
For saturated $u$, this is clear by the
construction of $\poldiva$.
For general degrees $u \in \gewkegel \cap M$, 
we have to argue as usual: If $g\in\K(Y)$, then
$$
g\in \Gamma(Y,\mathcal{O}_Y(\poldiva(u)))
\Leftrightarrow
g^n\in \Gamma(Y,\mathcal{O}_Y(\poldiva(nu)))
\Leftrightarrow
(gs(u))^n\in A_{nu}
\Leftrightarrow
gs(u)\in A_{u}.
$$
Note that for the last step one uses normality 
of the ring $A$.
\end{proof}

\section{Fibers of the quotient map}
\label{sec:fibers}

The construction of the affine $T$-variety
$X$ associated to a pp-divisor $\poldiva$ 
on a normal semiprojective variety $Y$  
involves, as an
intermediate step, the construction of a 
certain $T$-variety $\til{X}$ over $Y$.
The aim of this section is to describe 
the geometry of the fibers of the canonical 
map $\pi \colon \til{X} \to Y$ in terms
of the defining pp-divisor.

Fix a normal semiprojective 
variety $Y$; let $\gittera$ be a lattice, 
and $\kegela \subset \gittera_{\Q}$ 
a pointed cone.
As outlined in Section~\ref{section2},
any pp-divisor 
$\poldiva \in \PPDiv_{\Q}(Y,\kegela)$
gives rise to a sheaf of graded
$\mathcal{O}_{Y}$-algebras
\begin{eqnarray*}
\mathcal{A}
& = & 
\bigoplus_{u \in \kegela^\vee \cap \dualgittera}
\mathcal{O}(\poldiva(u)),
\end{eqnarray*}
where $M$ is the dual lattice of $N$.
The variety $\til{X} := \Spec_Y(\mathcal{A})$
comes with an action of the torus 
$T := \Spec(\K[M])$, and the canonical map 
$\pi \colon \til{X} \to Y$ is
a good quotient for the action of $T$
on $\til{X}$.

A first step in the study of the fibers 
of $\pi \colon \til{X} \to Y$ is an 
investigation of certain bouquets of 
toric varieties associated 
to $\sigma$-polyhedra 
$\polyedera \subset N_{\Q}$.
For the definition of these objects,
recall that one associates a normal quasifan
$\Lambda(\Delta)$ to $\polyedera$; the faces 
of $\polyedera$ correspond to the cones of 
$\Lambda(\Delta)$ via
\begin{eqnarray*}
F 
& \mapsto &
\lambda(F)
\; = \; 
\{u \in \dualgittera_{\Q}; \; 
\langle u, v - v' \rangle \ge 0
\text{ for all } 
v \in \polyedera, \, v' \in F\}.
\end{eqnarray*}
Having obtained a quasifan $\Lambda$ 
from the $\sigma$-polyhedron $\polyedera$,
we associate to this quasifan a graded algebra
$\K[\Lambda]$ by a frequently used procedure,
and then define the toric bouquet associated 
to $\polyedera$ as the 
spectrum of the graded algebra 
$\K[\Lambda]$.

\begin{definition}
Let $M$ be a lattice, and let $\Lambda$ be a
quasifan with convex support $\omega\subseteq M_\Q$. 
The {\em fan ring\/} associated
to $\Lambda$ is the affine $\K$-algebra 
defined by
$$
\K[\Lambda]
\; := \;
\bigoplus_{u \in \omega \cap M} \K \chi^{u},
\qquad
\chi^{u} \chi^{u'} :=
\left\{ 
\begin{array}{ll} 
\chi^{u+u'} & \text{ if } u,u' \in \lambda
              \text{ for some } \lambda \in  \Lambda, \\
0           & \text{ else }.
\end{array}
\right. 
$$
\end{definition}

Note that the fan ring $\K[\Lambda]$ may as well
be viewed as a semigroup algebra, if we define 
$u + u' := 0$, whenever $u$ and $u'$ do not belong
to a common cone.

\begin{definition}
Let $\Delta \subset N_\Q$ be a $\sigma$-polyhedron. 
The {\em toric bouquet\/} associated to $\Delta$ is
$X(\Delta) := \Spec(\K[\Lambda])$,
where $\Lambda=\Lambda(\Delta)$ is the normal quasifan of $\Delta$.
\end{definition}

We collect some basic geometric properties of
these toric bouquets; in particular, we 
note that they have equidimensional toric
varieties as their irreducible components, 
whence the name.
The proofs of the statements only use standard
toric geometry and therefore are left to the 
reader.

\begin{proposition}
\label{bouquets}
Let $\Delta \subset N_\Q$ be a $\sigma$-polyhedron,
let $\Lambda = \Lambda(\Delta)$ be the normal quasifan
of~$\Delta$, and let $X(\Delta) =  \Spec(\K[\Lambda])$
be the corresponding toric bouquet.
\begin{enumerate}
\item
The $M$-grading of $\K[\Lambda]$ defines an effective
algebraic action of the torus $T := \Spec(\K[M])$ on
$X(\Delta)$.
\item
The $T$-orbits of $X(\Delta)$ are in dimension
reversing one-to-one correspondence with the faces of
$\Delta$ via $F \mapsto T \mal x_{F}$,
where $x_{F} \in X(\Delta)$ is defined by
$$
\chi^{u}(x_F) 
\; = \;
\left\{ 
\begin{array}{ll} 
1           & \text{ if } u \in \lambda(F), \\
0           & \text{ else }.
\end{array}
\right. 
$$
\item
For a face $F \preceq \Delta$, let
$I(F) \subset \K[\Lambda]$ denote the ideal
generated by the $\chi^{u}$'s with
$u \not\in \lambda(F)$.
Then the closure of the orbit through $x_F$ is
given by
\begin{eqnarray*}
\blk{T \mal x_F} 
&  = &  
V(I(F)).
\end{eqnarray*}
This orbit closure is a toric
variety, and,
denoting by $\lin(F) \subset N_\Q$ the vector
space generated by all $v-v'$ with $v,v' \in F$,
its defining cone is
\begin{eqnarray*}
\Q_ {\ge 0} \mal (\Delta-F) 
/
\lin(F)
& \subset &
(N / (\lin(F) \cap N))_{\Q}. 
\end{eqnarray*}
\item
The irreducible components of $X(\Delta)$
are precisely the orbit closures corresponding to
the vertices $v \in \Delta$. As toric varieties,
they correspond to the cones
\begin{eqnarray*}
\lambda(v)^\vee 
& = &
\Q_{\ge 0} \mal (\Delta-v).
\end{eqnarray*}
\end{enumerate}
\end{proposition}

For our description of the fibers of the 
canonical map $\pi \colon \til{X} \to Y$, 
we have to generalize the toric bouquets 
discussed so far. 
We have to consider certain subalgebras 
of the fan ring arising from a 
$\sigma$-polyhedron.

\begin{definition}
\label{finitedef}
Let $\Lambda$ be the normal quasifan 
in $M_\Q$ 
of a $\sigma$-polyhedron 
$\Delta \subset N_\Q$. 
Let $\omega=\sigma^\vee$ and $S \subset \omega \cap M$ be a subset 
such that for each cone 
$\lambda \in \Lambda$ we have 
$$
S_{\lambda} 
\; := \; 
\lambda \cap S 
\; = \; 
\lambda \cap M_{\lambda}
$$
with a lattice 
$M_{\lambda} \subset M \cap \lin(\lambda)$ 
of full rank in $\lin(\lambda)$. 
Then we obtain a finitely generated subalgebra: 
$$
\K[\Lambda,S] 
\; := \;  
\bigoplus_{u \in S} \K \chi^{u} 
\; \subset \; 
\K[\Lambda].
$$
\end{definition}

\begin{example}
\label{varlattice1}
Let $N := \Z$, and $\sigma := \{0\}$.
Then any interval $\Delta \subset \Q$ 
is a $\sigma$-polyhedron, and the 
corresponding normal fan
$\Lambda$ consists of 
$$ 
\lambda_- := (-\infty,0], 
\qquad 
\lambda_0 := \{0\},
\qquad
\lambda_+ := 
[0, \infty).
$$
We have $\omega \cap M = \Z$, and
fixing a subset $S \subset \Z$ means to 
establish a ``conewise varying 
lattice structure''.
For example, we may obtain a picture

\medskip

\begin{center}
\begin{picture}(0,0)%
\includegraphics{omcomplex.pstex}%
\end{picture}%
\setlength{\unitlength}{1243sp}%
\begingroup\makeatletter\ifx\SetFigFont\undefined%
\gdef\SetFigFont#1#2#3#4#5{%
  \reset@font\fontsize{#1}{#2pt}%
  \fontfamily{#3}\fontseries{#4}\fontshape{#5}%
  \selectfont}%
\fi\endgroup%
\begin{picture}(7720,797)(416,-785)
\put(3511,-691){\makebox(0,0)[lb]{\smash{\SetFigFont{8}{9.6}{\familydefault}{\mddefault}{\updefault}{\color[rgb]{0,0,0}$0$}%
}}}
\put(2251,-691){\makebox(0,0)[lb]{\smash{\SetFigFont{8}{9.6}{\familydefault}{\mddefault}{\updefault}{\color[rgb]{0,0,0}$-2$}%
}}}
\put(4861,-691){\makebox(0,0)[lb]{\smash{\SetFigFont{8}{9.6}{\familydefault}{\mddefault}{\updefault}{\color[rgb]{0,0,0}$3$}%
}}}
\end{picture}

\end{center}

\noindent
where the subset $S \subset \Z$ is 
indicated by the fat points; it is  
is obtained by means of choice 
of sublattices 
$M_\lambda \subset \Z$ conpatible with
the fan structure: 
$$ 
S \cap \lambda_- := (-\infty,0] \cap 2\Z, 
\qquad 
S \cap \lambda_0 := \{0\},
\qquad
S \cap \lambda_+ := 
[0, \infty) \cap 3 \Z.
$$
\end{example}

Geometrically, the passage from $\K[\Lambda]$
to a subalgebra $\K[\Lambda,S]$ as in 
Definition~\ref{finitedef}
corresponds to taking componentwise quotients 
of a toric bouquet by a possibly varying 
finite group. More precisely, we obtain:

\begin{proposition}
\label{bouquets2}
Let $\Delta$, $\Lambda$, $S$, etc. be as in
Definition~\ref{finitedef}.
Then the inclusion
$\K[\Lambda,S] \subset \K[\Lambda]$
gives rise to a finite equivariant morphism
$\nu \colon X(\Delta) \to X(\Delta,S)$
of affine $T$-varieties. Moreover, the following
statements hold:
\begin{enumerate}
\item
$\blk{T \mal x_F} = \nu^{-1}\bigl(\blk{T \mal \nu(x_F)}\bigr)$
for every orbit closure $\blk{T \mal x_F} \subset X(\Delta)$,
and on $\blk{T \mal x_F}$, the map $\nu$ is the quotient by a
finite group $H_F$:
$$
\xymatrix{
{\blk{T \mal x_F}} \ar[rr]^{\nu} \ar[dr]
& & 
{\blk{T \mal \nu(x_F)}}
\\
& {\blk{T \mal x_F}}/H_F \ar[ur]_{\cong} & 
}
$$
Setting
$M_{F,S} 
 := 
(\lin(\lambda(F)) \cap M) 
/ 
M_{\lambda(F)}$,
the finite group $H_F$ is obtained as follows
$$
H_F 
\; = \; 
\Spec(\K[M_{F,S}]) 
\; = \; 
T_{\nu(x_{F})}/T_{x_F}.
$$
\item
The orbit closure $\blk{T \mal \nu(x_F)}$
is a toric variety with big torus $T/T_{\nu(x_F)}$.
Its corresponding cone is the image of
$\Q_ {\ge 0} \mal (\Delta-F)$ under
$N \to \Hom(M_{\lambda(F)},\Z)$.
\end{enumerate}
\end{proposition}

\begin{proof}
For every face $F \preceq \Delta$, we have the 
ideal $I(F,S) := \nu^*(I(F))$ in $\K[\Lambda,S]$.
The corresponding homomorphism of factor algebras
is the inclusion of a certain Veronese subalgebra:
\begin{eqnarray*}
\bigoplus_{u \in \lambda(F) \cap M_{\lambda(F)}} \K \chi^u
& \subset &
\bigoplus_{u \in \lambda(F) \cap M} \K \chi^u.
\end{eqnarray*}
In particular, the left hand side is the ring of
invariants of the group $H_F$. 
Moreover, by Proposition~\ref{bouquets}~(iii), 
this inclusion of algebras describes the map
$\nu$ over $\nu(\blk{T \mal x_F})$. 
The assertions follow directly. 
\end{proof}

We now associate a couple of fiber
data to any $y \in Y$. 
We say that a divisor $D \in \WDiv_\Q(Y)$
is {\em principal at $y \in Y$\/}
if near $y \in Y$ it is the divisor of a
function $f \in \K(Y)$.

\begin{definition}
\label{fibredatadef}
Let $\mathfrak{D} = \sum \Delta_i \otimes D_i$
be a representation such that all $D_i$ are prime,
and let $y \in Y$.
\begin{enumerate}
\item 
We define the  
{\em fiber polyhedron\/} 
of $y \in Y$ to be
$$
\Delta_y 
\;  :=  \;
\sum_{y \in D_i} \Delta_i
\; \in \;  
\Pol_{\sigma}(N_\Q).
$$
\item 
The normal quasifan of
$\Delta_y$ is denoted by 
$\Lambda_y$. 
\item
We define the
{\em fiber monoid complex\/} 
of $y \in Y$ as
\begin{eqnarray*}
S_y
& := &
\{
  u \in \omega \cap M; \; 
  \mathfrak{D}(u) \text{ is principal at } y
\}.
\end{eqnarray*}
\item
For $\lambda \in \Lambda_y$,
we denote by $M_{y,\lambda} \subset M$ the
sublattice generated by $S_y \cap \lambda$.
\end{enumerate}
\end{definition}

\begin{example}
\label{varlattice2}
Let $Y := \K^1$; 
take $N := \Z$,
and let $\sigma \subset \Q$
be the zero cone.
Consider the pp-divisor 
\begin{eqnarray*}
\poldiva
& = & 
\left[\frac{1}{3}, \frac{1}{2}\right] \otimes \{0\}.
\end{eqnarray*} 
Then, for $y = 0$, the fiber polyhedron
is $[1/3,1/2]$, and the fiber monoid 
complex is as the picture in Example~\ref{varlattice1}.
%
%
%
For any other point $y \in Y$, the fiber polyhedron
is $\sigma = \{0\}$, and the fiber monoid complex
is just $\Z$.
\end{example}

\begin{lemma}
\label{gordon}
Let $y \in Y$. 
Then every $\lambda \in \Lambda_y$
satisfies 
$S_y \cap \lambda = M_{y,\lambda} \cap \lambda$. 
In particular, $S_y \cap \lambda$ is a finitely 
generated semigroup, and 
we obtain a finitely generated algebra
$$
\K[\Lambda_y,S_y]
\; := \; 
\bigoplus_{u \in S_y} \K \chi^u 
\; \subset \; 
\K[\Lambda_y].
$$
\end{lemma}

Now we are ready to begin with the 
study  
of the fibers of the canonical map 
$\pi \colon \til{X} \to Y$.
Below is the first statement.

\begin{proposition}
\label{pifibers}
Let $\mathfrak{D} \in \PPDiv_{\Q}(Y,\sigma)$.
Then, for every $y \in Y$, the reduced fiber 
$\pi^{-1}(y)$ of the associated map 
$\pi \colon \til{X} \to Y$ is 
$T$-equivariantly isomorphic to 
$X(\Delta_y,S_y) = \Spec(\K[\Lambda_y,S_y])$ 
where $\Lambda_y$ is the normal quasifan of 
$\Delta_y$. 
\end{proposition}

\begin{proof}
Fix $y \in Y$, and consider the group 
$G \subset \CDiv_{\Q}(Y)$ generated by 
those $\mathfrak{D}(u)$
with $u \in \omega \cap M$ that are 
principal
at $y$. 
By Lemma~\ref{gordon}, $G$ is finitely generated. 
Thus, we may choose a basis 
$E_1, \ldots, E_r \in G$. 
After possibly shrinking $Y$, we may assume that
$Y$ is affine, and that $E_i = \kdiv(g_i)$ holds
with $g_i \in \K(Y)$. 
For $D \in G$, set 
$$
g_D
\; :=  \;
g_1^{a_1} \ldots g_r^{a_r}, 
\qquad
\text{where }
D = a_1E_1 + \ldots + a_rE_r. 
$$

These functions satisfy
$g_{D'+D} = g_{D'}g_{D}$ for 
all $D',D \in G$. 
Consequently, denoting by $\Lambda_y$ 
the normal quasifan of $\Delta_y$,
we may define a graded epimorphism
$$
\Phi 
\colon
\Gamma(Y,\mathcal{A}) 
\to 
\K[\Lambda_y,S_y],
\qquad
\Gamma(Y,\mathcal{A})_u
\ni h
\mapsto 
\left\{ 
\begin{array}{ll} 
(g_{\mathfrak{D}(u)}h)(y) \chi^{u}           
& \text{if } u \in S_y, \\
0           
& \text{else}.
\end{array}
\right. 
$$
To see multiplicativity, let
$h_i \in \Gamma(Y,\mathcal{A})_{u_i}$.
Then, $g_{\mathfrak{D}(u_1+u_2)}h_1h_2$
vanishes at $y$ if the $u_i$ 
do not belong to the same $\lambda \in \Lambda$, 
or if one 
of the $\mathfrak{D}(u_i)$ is not principal
at $y$. 

To conclude the proof, we have to show
that the kernel of the above 
epimorphism $\Phi$ equals the radical of 
the ideal of the fiber $\pi^{-1}(y)$. 
The fiber ideal is given by  
$$
I_y
\; := \;
\langle 
hf; 
\; 
h \in \Gamma(Y,\mathcal{O}),
\; 
h(y) = 0,
\;
f \in \Gamma(Y,\mathcal{A})_u, 
\; 
u \in \omega \cap M
\rangle
\; \subset \; 
\Gamma(Y,\mathcal{A}). 
$$
Obviously, $I_y \subset \ker(\Phi)$. 
Conversely, for any $f \in \ker(\Phi)$, 
say homogeneous of degree $u$, 
we have $f^n \in I_y$ as 
soon as $\mathfrak{D}(nu)$ 
is principal at $y$.
\end{proof}

Combining the above result with
Propositions~\ref{bouquets} and~\ref{bouquets2} 
and the notions of Definition \ref{orbitdatadef}
gives the following information on the 
geometry of the fibers of  
$\pi \colon \til{X} \to Y$:

\begin{corollary}
\label{orbitsfiber}
For $y \in Y$, consider the affine 
$T$-variety $\pi^{-1}(y)$.
\begin{enumerate}
\item  
For any $\til{x} \in \pi^{-1}(y)$, 
we have $\omega(\til{x}) \in \Lambda_y$,
and this sets up a one-to-one correspondence 
between the $T$-orbits of $\pi^{-1}(y)$ 
and the cones of $\Lambda_y$
(corresponding to the faces of $\Delta_y$).
The orbit lattice of $\til{x}$ 
is $M(\til{x}) = M_{y,\omega(\til{x})}$. 
\item 
The irreducible components of 
$\pi^{-1}(y)$ are the orbit
closures $\blk{T \mal \til{x}}$
with $\omega(\til{x})$ maximal in
$\Lambda_y$. They are 
normal toric varieties with 
big torus $T/T_{\til{x}}$, where 
$T_{\til{x}} = \Spec(\K[M/M_{y,\omega(\til{x})}])$.
\end{enumerate}
\end{corollary}

In the proof of Proposition~\ref{pifibers},
we had to compare the fiber ideal $I_y$ and its
radical $\ker(\Phi)$. 
Looking a little bit closer at these
data gives the following:

\begin{proposition}
\label{reduced}
A fiber $\pi^{-1}(y)$ is reduced if and only if 
all $\mathfrak{D}(u)$, where $u \in \omega \cap M$
are principal at $y \in Y$. 
\end{proposition} 

\begin{example}
As in Example~\ref{varlattice2}, 
take $Y := \K^1$,
let $N := \Z$,
and $\sigma := \{0\}$.
Again consider the pp-divisor 
\begin{eqnarray*}
\poldiva
& = & 
\left[\frac{1}{3}, \frac{1}{2}\right] \otimes \{0\}.
\end{eqnarray*} 
Then the $\K^*$-variety $\til{X}$ associated
to $\poldiva$ is the affine space $\K^2$
together with the $\K^*$-action given
by
\begin{eqnarray*}
t \mal (z,w)
& = & 
(t^3z,t^{-2}w).
\end{eqnarray*}
The canonical map $\pi \colon \til{X} \to Y$ 
may, as a good quotient for this $\K^*$-action,
concretely be written as
$$ 
\pi \colon
\K^2 \; \mapsto \; \K,
\qquad
(z,w) \; \mapsto \; z^{2}w^3.
$$     
The fiber $\pi^{-1}(0)$ is the union of three 
orbits: the origin and the orbits through
$(1,0)$ and $(0,1)$.

Combinatorially, this is reflected as follows.
We have
$\Delta_0 = [1/3, 1/2]$ as the fiber polyhedron
in $0 \in \K$;
the associated normal fan $\Lambda_0$ consists 
of three cones as in Example~\ref{varlattice1}, 
and the fiber monoid complex $S_0$ has a 
varying lattice structure:

\medskip

\begin{center}
\begin{picture}(0,0)%
\includegraphics{omcomplex.pstex}%
\end{picture}%
\setlength{\unitlength}{1243sp}%
\begingroup\makeatletter\ifx\SetFigFont\undefined%
\gdef\SetFigFont#1#2#3#4#5{%
  \reset@font\fontsize{#1}{#2pt}%
  \fontfamily{#3}\fontseries{#4}\fontshape{#5}%
  \selectfont}%
\fi\endgroup%
\begin{picture}(7720,797)(416,-785)
\put(3511,-691){\makebox(0,0)[lb]{\smash{\SetFigFont{8}{9.6}{\familydefault}{\mddefault}{\updefault}{\color[rgb]{0,0,0}$0$}%
}}}
\put(2251,-691){\makebox(0,0)[lb]{\smash{\SetFigFont{8}{9.6}{\familydefault}{\mddefault}{\updefault}{\color[rgb]{0,0,0}$-2$}%
}}}
\put(4861,-691){\makebox(0,0)[lb]{\smash{\SetFigFont{8}{9.6}{\familydefault}{\mddefault}{\updefault}{\color[rgb]{0,0,0}$3$}%
}}}
\end{picture}

\end{center}

\noindent
This varying lattice structure indicates that 
$\K^*$ acts on one
coordinate axis of $\K^2$ with generic isotropy group 
of order two,
and on the other with generic isotropy group of 
order three.
\end{example}

\section{Functoriality properties}

Our first results, Theorems~\ref{data2variety} 
and~\ref{variety2data}, establish correspondences between
pp-divisors on semiprojective varieties,
on the one hand, and affine varieties with effective 
torus action on the other.
In this section, we present the functoriality properties
of these assignments;
the proofs of the results are given in 
Section~\ref{mainamainbpf}. 

Going from polyhedral divisors to varieties is 
functorial in an almost evident manner, 
but the reverse direction is more delicate.
Nevertheless, in an appropriate setup,
we obtain an equivalence of categories,
and our results allow us to decide 
when two pp-divisors define isomorphic
$T$-varieties,
see Corollaries~\ref{equivcat} and~\ref{isoclasses}.

First, we have to fix the respective notions of 
morphisms. Concerning varieties with torus action, 
we will work with the following concept.

\begin{definition}
\label{eqmorphdef}
Let $X$ and $X'$ be varieties endowed with effective 
actions of tori $T$ and $T'$.
By an {\em equivariant morphism\/} 
from $X$ to $X'$, from now on we mean a morphism 
$\morphX \colon X  \to  X'$ admitting an 
accompanying homomorphism 
$\til{\morphX} \colon T \to T'$ 
such that  $\morphX(t \mal x) = \til{\morphX}(t) \mal \morphX(x)$ 
holds  for all $(t,x) \in T \times X$.
\end{definition}

So, a morphism of two $T$-varieties 
is equivariant in the usual 
sense, if and only if it has the 
identity as an accompanying homomorphism.
Note that in case of a dominant
morphism $\morphX \colon X  \to  X'$,
the accompanying homomorphism is 
uniquely determined.

We turn to pp-divisors. To define
the notion of a map between two pp-divisors,
we first have to introduce the concept
of a ``polyhedral principal divisor''.

\begin{definition}
\label{polfuncdef}
Let $Y$ be a normal variety,
$\gittera$ a lattice and  
$\kegela \subset \gittera_{\Q}$
a pointed cone.
\begin{enumerate}
\item
A {\em plurifunction\/} 
with respect to the lattice $\gittera$
is an element of 
$\K(Y,N)^* := \gittera \otimes_\Z \K(Y)^*$.
\item 
For $u \in \dualgittera = \Hom(\gittera,\Z)$,
the {\em evaluation\/} of a plurifunction 
$\polfunca = \sum v_i \otimes f_i$ 
with respect to $\gittera$ is
$$
\polfunca(u)
\; := \; 
\prod f_i^{\langle u,v_i \rangle}
\; \in \;
\K(Y)^*.
$$
\item
The {\em ``polyhedral principal'' divisor\/} 
with respect to $\kegela \subset \gittera_\Q$ 
of a plurifunction 
$\polfunca = \sum v_i \otimes f_i$ 
with respect to $\gittera$ is
$$
\ddiv(\polfunca) 
\; := \;
\sum (v_i + \kegela) \otimes \ddiv(f_i)  
\; \in \;
\CDiv(Y,\kegela).
$$
\end{enumerate}
\end{definition}

Note that evaluating and taking the 
divisor with respect to $\kegela \subset \gittera_\Q$
of a plurifunction commute for 
$u \in \kegela^\vee \cap \dualgittera$.
We are ready to define the notion 
of a map of pp-divisors.

\begin{definition}
\label{homprp}
Let $Y$, $Y'$ be normal semiprojective varieties, 
$\gittera$, $\gittera'$ lattices, 
$\kegela \subset \gittera_\Q$ and
$\kegela' \subset \gittera'_\Q$
pointed cones, 
and consider pp-divisors
$$
\poldiva 
 \; = \;
\sum \Delta_i \otimes D_i  
\; \in \; 
\PPDiv_{\Q}(Y, \kegela),
\qquad
\poldiva' 
\; = \;
\sum \Delta'_i \otimes D'_i  
\; \in \;
\PPDiv_{\Q}(Y',\kegela').
$$ 
\begin{enumerate}
\item
For morphisms $\morphY \colon Y \to Y'$ such that
none of the supports $\Supp(D_i')$
contains $\morphY(Y)$, we define the 
(not necessarily proper)
{\em polyhedral pull back\/} as
$$
\morphY^*(\poldiva') 
\; := \;
\sum \polyedera_i' \otimes  \morphY^*(D_i') 
\; \in \;
\CDiv_\Q(Y,\kegela').
$$
\item
For linear maps $\linabba \colon \gittera \to \gittera'$ 
with $\linabba(\kegela) \subset \kegela'$,
we define the (not necessarily proper) 
{\em polyhedral push forward\/} as
$$
\linabba_*(\poldiva) 
\; := \; 
\sum (\linabba (\polyedera_i) + \kegela') \otimes D_i 
\; \in \;
\WDiv_\Q(Y,\kegela').
$$
\item
A {\em map\/} $\poldiva \to \poldiva'$ is a triple 
$(\morphY,\linabba,\polfunca)$ with a dominant
morphism $\morphY \colon Y \to Y'$,
a linear map 
$\linabba \colon \gittera \to \gittera'$
as in~(ii), and a plurifunction
$\polfunca \in \K(Y,\gittera')^{*}$ 
such that
\begin{eqnarray*}
\morphY^*(\poldiva') 
& \le & 
\linabba_*(\poldiva) 
+ \ddiv(\polfunca).
\end{eqnarray*}
\end{enumerate}
\end{definition}

Note that the relation ``$\leq$'' among pp-divisors
is equivalent to requiring 
the opposite inclusion for the respective 
polyhedral coefficients 
in the representation of the pp-divisors
as linear combinations of prime divisors.

\begin{example}
\label{adjunctionmaps}
With the notation of the previous definition,
we obtain natural adjunction maps:
\begin{enumerate}
\item 
For any generically finite morphism
$\morphY \colon Y \to Y'$, 
the pullback $\psi^{*}(\poldiva')$
is, if defined at all, a pp-divisor on $Y$, and the triple
$(\psi,\id_{\gittera},1)$ defines a map 
$\psi^{*}(\poldiva') \to \poldiva'$.
\item For any lattice homomorphism 
$\linabba \colon \gittera \to \gittera'$, 
the triple $(\id_{Y},\linabba,1)$ defines a map
$\poldiva \to \linabba_{*}(\poldiva)$,
provided that $\linabba_{*}(\poldiva)$ is a 
pp-divisor.
\end{enumerate}
\end{example}

In order to obtain a category of pp-divisors,
we still have to introduce composition.
For this, note that along the lines 
of~\ref{homprp}~(i) and~(ii),
we can also define pullback and pushforward of 
plurifunctions.

\begin{definition} 
The {\em identity map\/} 
of a pp-divisor is the triple 
$(\id,\id,1)$. 
The {\em composition\/} of two maps of pp-divisors
$(\morphY,\linabba,\polfunca)$ and $(\morphY',\linabba',\polfunca')$
is defined as the map of pp-divisors
$(\morphY' \circ \morphY, \, \linabba' \circ \linabba, \,
F'_{*}(\polfunca) \cdot \psi^{*}(\polfunca'))$.
\end{definition}

Let us now demonstrate how to obtain functoriality.
The construction of Theorem~\ref{data2variety} associates to
a given pp-divisor $\poldiva \in \PPDiv_{\Q}(Y,\kegela)$ 
on a normal semiprojective variety $Y$ the affine variety
$$
\funcX(\poldiva)
\; := \; 
X \; := \;
\Spec (\Gamma(Y,\mathcal{A})),
\quad
\text{where }
\mathcal{A}
\; = \; 
\bigoplus_{u \in \sigma^\vee \cap M} 
\mathcal{O}(\poldiva(u)).
$$

\begin{proposition}
\label{prop:functorial}
The assignment $\poldiva \mapsto \funcX(\poldiva)$ is a 
faithful covariant functor from the pp-divisors 
on normal semiprojective varieties
to the normal affine varieties with torus action. 
\end{proposition}

\begin{proof}
Let $\poldiva$ and $\poldiva'$ be pp-divisors on  
$Y$ and $Y'$ respectively. 
Write $X := \funcX(\poldiva)$ and $X' := \funcX(\poldiva')$,
and denote the acting tori by $T$ and $T'$.
By definition, any map 
$\poldiva \to \poldiva'$,
given by $(\morphY,\linabba,\polfunca)$,
induces homomorphisms of 
$\Gamma(Y',\mathcal{O})$-modules:
\begin{equation}
\label{eq:functor}
\Gamma(Y',\mathcal{O}(\poldiva'(u)))
\; \to \;
\Gamma(Y,\mathcal{O}(\poldiva(\linabba^{*} u))),
\qquad
h 
\; \mapsto \;
\polfunca(u) \morphY^*(h).
\end{equation}
These maps fit together to a graded homomorphism
$\Gamma(Y',\mathcal{A}') \to \Gamma(Y,\mathcal{A})$.
This in turn gives an equivariant morphism 
$\morphX \colon X \to X'$
with the map 
$\til{\morphX} \colon T \to T'$ defined by 
$\linabba \colon \gittera \to \gittera'$
as an accompanying homomorphism.

Obviously, the identity map of pp-divisors
defines the identity on the level of
equivariant morphisms. 
Compatibility with composition follows
from the definition of the 
equivariant morphism associated to 
a map of pp-divisors via~\ref{eq:functor}
and the fact that we always have
\begin{eqnarray*}
F'_{*}(\polfunca) \cdot \psi^{*}(\polfunca')(u)
& = & 
\polfunca((F')^*(u)) \psi^{*}(\polfunca'(u)).
\end{eqnarray*}

In order to see that the functor is faithful, 
i.e.\ injective on 
morphisms, consider two maps  
$(\morphY_i,\linabba_i,\polfunca_i)$
of pp-divisors $\poldiva$ on $Y$ 
and $\poldiva'$ on $Y'$ that define the 
same equivariant morphism 
$\morphX \colon X \to X'$.

To obtain $\linabba_1 = \linabba_2$,
it suffices to check that $\varphi(X)$ contains points 
with free $T'$-orbit: by equivariance, 
this will fix the accompanying homomorphism
$\til{\varphi} \colon T \to T'$, which in turn determines 
the lattice homomorphisms $F_i \colon N \to N'$.
 
According to the properties of the maps $\morphY_i$, 
there is an open 
set $V \subset Y'$ with $\morphY_i(Y) \cap V \ne \emptyset$
such that any $u \in \omega \cap \dualgittera$ 
admits a section $f_u \in \mathcal{A}_u(V)$ without
zeroes in $V$. 
Consider
\begin{eqnarray*}
U
& := &
\bigcap X \setminus \ddiv(\varphi^*(f_u)).
\end{eqnarray*}
Using Theorem~\ref{data2variety}~(iii), one sees that $U$
is nonempty. 
Moreover, for each $x \in U$ and each $u \in M$, the image 
$\varphi(x)$ admits a homogeneous rational function 
$f \in \K(X')_u$ defined near $\varphi(x)$ with $f(\varphi(x)) \ne 0$. 
In other words, $\varphi(x)$ has a free $T'$-orbit. 

In order to see $\morphY_1 =  \morphY_2$,
it suffices to show that 
$\morphY_1^*, \morphY_2^* \colon \K(Y') \to \K(Y)$
coincide.
Given a function $f \in \K(Y')$, we may write it
as $f = g/h$ with $g,h \in \Gamma(Y',\mathcal{A}_u)$
for some $u \in \dualgittera'$, because we have 
$\K(Y') = \K(X')^T$, and on $X'$ any invariant
rational function is the quotient of two semiinvariants.
Thus, we obtain
$$
\morphY_i^*(f)
\ = \ 
\frac{\polfunca_i(u) \morphY_i^*(g)}{\polfunca_i(u) \morphY_i^*(h)}
\ = \ 
\frac{\varphi^*(g)}{\morphX^*(h)}.
$$

Finally, the equality $\polfunca_1 = \polfunca_2$
follows from their appearance in the comorphism 
of  $\morphX \colon X \to X'$:
since the weight monoid of the $\torusa'$-action 
on $X'$ generates the lattice $\dualgittera'$,
there are enough homogeneous sections $h$
to fix $\polfunca_i$ by using the defining
formula~\ref{eq:functor} of $\varphi^*$ and 
$\psi_1=\psi_2$. 
\end{proof}

To proceed, we introduce a notion of minimality for 
a pp-divisor $\poldiva \in \PPDiv_\Q(Y,\sigma)$.
Since all evaluations $\poldiva(u)$ are 
semiample, they come along with a natural
contraction map, being birational whenever
$u \in \relint{\sigma^\vee}$:
$$
\vartheta_u \colon Y 
\; \to \; 
Y_u
:= 
\Proj\left(
\bigoplus_{n \in \Z_{\ge 0}} 
\Gamma(Y,\mathcal{O}(\poldiva(nu)))
\right).
$$ 
Denoting by $X$ the normal, affine $T$-variety
associated to the pp-divisor $\mathfrak{D}$, 
we can recover the semiprojective varieties $Y_u$ 
as GIT-quotient spaces associated to linearizations 
of the trivial bundle; namely, we have
\begin{eqnarray*}
Y_u 
& = & 
X^{ss}(u) \quot T.
\end{eqnarray*}
{From} this we see, similarly as in Section~\ref{limit},
that the spaces $Y_u$ fit into an inverse system
with projective morphisms 
$\vartheta_{uw} \colon Y_u \to Y_w$, whenever 
$u \in \lambda$ and $w \in \gamma$ for two 
cones $\gamma \preceq \lambda$ of the GIT-quasifan
of the $T$-action on $X$.  
  
Clearly, we have $\vartheta_w = \vartheta_{uw} \circ \vartheta_u$,
whenever composition is possible. Thus, the 
morphisms
$\vartheta_u \colon Y \to Y_u$ lift to a 
(projective, birational) 
morphism to the inverse limit of the system of the 
GIT-quotient spaces:
$$ 
\vartheta \colon Y \; \to \; \kprojlim Y_u.
$$
Recall from Section~\ref{limit} that $\kprojlim Y_u$ 
comes with a canonical component, dominated by the 
intersection $W$ over all $X^{ss}(u)$, where 
$u \in \sigma^\vee \cap M$.
By construction, $\vartheta$ maps $Y$ 
onto this component.

\begin{definition}
We say that a pp-divisor 
$\poldiva \in \PPDiv_\Q(Y,\sigma)$
is {\em minimal\/} if the morphism
$\vartheta \colon Y \to \kprojlim Y_u$ is the 
normalization of the canonical component of 
$\kprojlim Y_u$.
\end{definition}

Note that the pp-divisors constructed in the
proof of Theorem~\ref{variety2data} are minimal.
Moreover, on a curve $Y$, every pp-divisor 
is minimal.

The following result makes precise, up to what 
extent we can describe equivariant morphisms 
in terms of maps of pp-divisors.

\begin{theorem}
\label{mainab}
Let
$
\,\poldiva 
\; \in \; 
\PPDiv_\Q(Y,\sigma)$ and
$\,\poldiva'
\; \in \; 
\PPDiv_\Q(Y',\sigma')
$
be pp-di\-vi\-sors
and let  
$\varphi \colon \funcX(\poldiva) \to \funcX(\poldiva')$
be a dominant,
equivariant morphism.
Then, there exist a projective birational
morphism $\kappa \colon \til{Y} \to Y$,
a map $(\psi,F,\polfunca)$ from
$\kappa^*{\poldiva}$ to $\poldiva'$,
and a commutative diagram
$$ 
\xymatrix{
& 
{\funcX(\kappa^*{\poldiva})} 
\ar[dl]_{\funcX(\kappa,\id,1)}^{\cong}  
\ar[dr]^{\funcX(\psi,F,\polfunca)} 
& 
\\ 
{\funcX(\poldiva)}  
\ar[rr]_{\varphi} 
& &
{\funcX(\poldiva')}. 
}
$$
If $\varphi$ is an isomorphism and
$\poldiva'$ is minimal,
then one may take $\kappa$ as the identity
and obtains $\varphi=\funcX(\psi,F,\polfunca)$
where $F \colon N \to N'$ is an
isomorphism of lattices sending $\sigma$
to $\sigma'$, and $\psi \colon Y \to Y'$ is
birational and projective; if also 
$\poldiva$ is minimal, then  
$\psi$ is an isomorphism.
\end{theorem}

The theorem shows that minimal 
pp-divisors may serve as a tool for the study of
equivariant automorphism groups,
where equivariant is understood in the usual
sense, i.e., with $\til{\varphi}=\id$ 
in the language of Definition~\ref{eqmorphdef}.

\begin{corollary}
Let $\poldiva \in \PPDiv_\Q(Y,\sigma)$ 
be a minimal pp-divisor.
Then the automorphisms 
$\varphi \colon \funcX(\poldiva) \to \funcX(\poldiva)$
satisfying $\varphi(t \mal x) = t \mal \varphi(x)$
correspond to pairs $(\psi,\polfunca)$,
where  $\psi \colon Y \to Y$ is an automorphism
and $\polfunca \in \K(N;Y)^*$ satisifies
$\psi^*(\poldiva) = \poldiva + {\rm div}(\polfunca)$.
\end{corollary}

\begin{example}
Given a lattice $N$ and a cone 
$\sigma \subset N_\Q$, we obtain
the associated affine toric variety
as $X_\sigma = \funcX(\poldiva)$ for the 
trivial pp-divisor $\poldiva=0$ 
living on $Y = \{y\}$,
see Example~\ref{afftorvar}.
For the big torus $T \subset X_\sigma$,
we have canonical identifications
$$ 
T 
\; \cong \; 
N\otimes_{\Z}\K^\ast
\; \cong \; 
\K(N;Y)^*.
$$
Hence, the translation $X_\sigma \to X_\sigma$
by a torus element $t \in T$ is the equivariant 
morphism $\funcX(\poldiva) \to \funcX(\poldiva)$
associated to the map 
$(\id_Y,\id_N,\polfunca)$
of $\poldiva$, 
where $\polfunca \in \K(N;Y)^*$ is the 
plurifunction corresponding to  $t \in T$.
\end{example}

\begin{example}
Let $\poldiva=\sum_{y\in Y} \Delta_y \otimes \{y\}$
be a pp-divisor on an elliptic curve $Y$.
In order to obtain the equivariant automorphisms
of $\funcX(\poldiva)$, we have to figure
the automorphisms $\psi \colon Y \to Y$
satisfying
$$
\psi^*(\poldiva) - \poldiva
\; = \; 
\sum_{y \in Y} \left(\Delta_{\psi(y)}- \Delta_y\right) \otimes \{y\}
\; = \;
\ddiv({\polfunca})
$$
with a plurifunction $\polfunca \in \K(N;Y)^*$.
Note that by completeness of $Y$,
the plurifunction $\polfunca \in \K(N;Y)^*$
is determined by its divisor up to a 
 ``constant'' from $N\otimes_{\Z}\K^\ast \cong T$;
so, $\ddiv(\polfunca)$ determines the automorphism
up to translation by a torus element.

The left thand side difference is a polyhedral 
principal divisor if and only if there are
elements $v_y\in N$ such that
$\Delta_{\psi(y)}=\Delta_y + v_y$ as 
polyhedra,
$\sum v_y=0$ in $N$, and,
using the group law on the elliptic curve,
$\sum v_y \otimes y =0$ in $N \otimes_\Z Y$.
In particular, unless 
$\Delta_y \in N + \sigma$ for all $y \in Y$, 
the automorphism $\psi$ must be of finite order.
\end{example}


As another immediate consequence, 
we can answer the question, when two given pp-divisors
define equivariantly isomorphic varieties.

\begin{corollary}
Two pp-divisors
$\poldiva_i \in \PPDiv_\Q(Y_i,\sigma_i)$
define equivariantly isomorphic varieties
$\funcX(\poldiva_i)$ if and only if there
are projective birational
morphisms $\psi_i \colon Y_i \to Y$ and
a pp-divisor
$\poldiva \in \PPDiv_\Q(Y,\sigma)$
with 
$\poldiva_i \cong \psi_i^* \poldiva$.
\end{corollary}

In order to turn the functor $\funcX$ into an equivalence of categories,
we restrict ourselves to those maps of polyhedral
divisors that define dominant equivariant morphisms.
Let us call these for the moment {\em dominating}.

\begin{remark}
If $(\morphY,\linabba,\polfunca)$ is a map of pp-divisors
$\poldiva$ and $\poldiva'$ such that $\linabba$ has finite
cokernel, then $(\morphY,\linabba,\polfunca)$ is dominating.
\end{remark}

However, the main obstruction for $\funcX$ to yield an equivalence
is the fact that, for a projective, birational map
$\psi:Y\to Y'$, the morphism $\psi^\ast(\poldiva')\to\poldiva'$
is not an isomorphism, but 
$\funcX(\psi^\ast(\poldiva'))\to\funcX(\poldiva')$ is.  
Hence, similar to the construction process of derived categories,
we have to localize by those maps: We extend the morphisms of our category of
pp-divisors by formally introducing an inverse of 
$\psi^\ast(\poldiva')\to\poldiva'$. The correct way to do this is to define
a new morphism $\poldiva'\to\poldiva''$ as a diagram 
\vspace{-1ex}
$$
\xymatrix{
\psi^\ast(\poldiva')\ar[d]\ar[dr]\\
{\poldiva'} & {\poldiva''}
}
$$
of traditional ones with some projective, birational $\psi$.
Now, as an immediate consequence of 
Theorem~\ref{mainab}, 
we obtain the following statement.

\begin{corollary}
\label{equivcat}
The functor $\funcX$ induces an equivalence from
the localized category of pp-divisors with dominating maps
to the category of normal affine varieties with
effective torus action and dominant equivariant morphisms.
\end{corollary}
 
Finally, we want to describe the isomorphism classes of normal,
affine $T$-varieties. Fixing $T$ and the weight cone $\omega$ of
its action, there are only two types of isomorphisms in
our localized category of pp-divisors.
First, the adding of divisors of plurifunctions: As we do
with principal divisors in the traditional setting, we handle this by
introducing the Picard group.

\begin{definition}
\label{defpPic}
Dividing by the group of polyhedral principal divisors,
we may define
the {\em polyhedral Picard group\/} and the
{\em rational polyhedral Picard group\/} as
$$
\Pic(Y,\kegela) 
 :=  
\CDiv(Y,\kegela) / \K(Y,N)^*,
\qquad
\Pic_{\Q}(Y,\kegela)
 := 
\CDiv_{\Q}(Y,\kegela) / \K(Y,N)^* .
$$
\end{definition}

Note that, by abuse of notation, we actually divide by the image
of $ \K(Y,N)^*$ -- there is always a kernel. Moreover, 
the rational polyhedral Picard group is {\em not\/}
the rational vector space associated to the (integral) polyhedral
Picard group. 

\begin{example}[Cf.~\ref{usual2polyhedral}]
Let $\gittera = \Z$ and $\kegela = \Q_{\ge 0}$.
Then we have $\Pol_{\kegela}(\gittera) = \Z$,
and the polyhedral Picard group is 
the usual one, i.e.\ $\Pic(Y,\kegela) = \Pic(Y)$.
Moreover, the following sequence is exact:
$$
0\to\Pic(Y)\to \Pic_\Q(Y,\kegela)\to \CDiv(Y)\otimes_\Z \Q/\Z\to 0.
$$
\end{example}

The second type of isomorphisms in the 
localized category of pp-divisors consists of the new isomorphisms
coming from (birational) modifications of $Y$. 
This yields 

\begin{corollary}
\label{isoclasses}
The isomorphism classes of normal affine varieties 
with effective $T$-action and fixed weight cone $\omega=\sigma^\vee$
are in 1-1 correspondence with the pp-classes
in $\kdirlim \Pic_\Q(\bullet,\kegela)$
where the limit is taken over modifications 
of the variety carrying the respective pp-divisors.
\end{corollary}

\begin{remark}
The normal affine $\Ta$-varieties $X$
with $\dim(\Ta) = \dim(X)-1$
are precisely those arising from
pp-divisors on smooth curves. 
Since there are no non-trivial modifications in the 
curve case, no localization is needed in this case.
\end{remark}

\section{Proof of Theorem~\ref{mainab}}
\label{mainamainbpf}

We begin with two auxiliary statements;
the first one is an elementary general 
observation on semiample divisors.

\begin{lemma}
\label{compdiv}
Let $D$ and $D'$ be semiample $\Q$-Cartier divisors on 
a normal variety~$Y$. 
If $\Gamma(Y,\mathcal{O}(nD)) \subset \Gamma(Y,\mathcal{O}(nD'))$ 
holds for infinitely many $n > 0$, then we have $D \leq D'$.
\end{lemma}

\begin{proof}
Write $D = \sum \alpha_iD_i$ and $D' = \sum \alpha_i'D_i'$
with prime divisors $D_i$ on $Y$.
Then, for any $y \in Y$, we obtain new divisors
by removing all prime components from $D$ and $D'$ 
that do not contain $y$:
$$
D_y \; := \; \sum_{y \in D_i} \alpha_iD_i,
\qquad 
D'_y \; := \; \sum_{y \in D'_i} \alpha'_iD'_i.
\qquad 
$$
By semiampleness of $D$, there are an $n > 0 $
and a global section $f \in \Gamma(Y,\mathcal{O}(nD))$
with $\ddiv(f)_{y} + nD_{y} = 0$. 
Since $f$ is a global section of 
$\mathcal{O}(D')$ as well, we have 
$\ddiv(f)_{y} + nD'_{y} \ge 0$.
This gives $D_{y} \le D'_{y}$ for
all $y \in Y$, which implies $D \le D'$. 
\end{proof}

Let us call the minimal 
pp-divisors as produced in the proof 
of Theorem~\ref{variety2data}
for the moment {\em GIT-constructed\/}.
So, given any $X$ (or a pp-divisor 
$\poldiva \in \PPDiv_\Q(Y,\sigma)$ defining
$X$), the associated GIT-constructed 
pp-divisors are minimal, and live 
on the normalization $\blk{Y}$ of the 
canonical component of the limit over 
the GIT-quotients in question.

\begin{lemma}
\label{GITconstrpb}
Let $\poldiva \in \PPDiv_\Q(Y,\sigma)$.
Then, for every associated GIT-constructed
$\blk{\poldiva} \in \PPDiv_\Q(\blk{Y},\sigma)$,
we have 
$\poldiva = \vartheta^* \blk{\poldiva} + {\rm div}(\polfunca)$ 
with a plurifunction $\polfunca \in \K(Y,N)^*$,
where $\vartheta \colon Y \to \blk{Y}$ is the canonical
morphism.
\end{lemma}

\begin{proof}
As usual, denote by 
$\mathcal{A}$ the $\mathcal{O}_{Y}$-algebra
associated to $\poldiva$,
let $\til{X} := \Spec_Y(\mathcal{A})$, 
set $A := \Gamma(Y,\mathcal{A})$,
and $X := \Spec(A)$.
Then, denoting by $r \colon \til{X} \to X$ 
the contraction map, we have 
commutative diagrams
\[
\label{catquot} 
\xymatrix{
{r^{-1}(X^{ss}(u))} \ar[r]^r \ar[d] & X^{ss}(u) \ar[d] 
\\
Y \ar[r]^{\vartheta_u} 
\ar[dr]_{\vartheta}
& Y_u
\\
&
{\blk{Y}} \ar[u]_{p_u}
}
\]

Now we are ready to compare $\poldiva$ and the 
pullback $\vartheta^{*}(\blk{\poldiva})$. 
For this, recall that the divisors 
$\blk{\poldiva}(u)$ have been defined  
in the proof of Theorem~\ref{variety2data}
via
$$
\mathcal{O}(\blk{\poldiva}(u))
\; := \; 
\frac{1}{s(u)} \cdot \blkk{\mathcal{A}}_u
\; \subset \;
\K(\blk{Y})
$$
with $s \colon \dualgittera \to \K(X)$ being a section 
of the degree ``map'',
and $\blkk{\mathcal{A}}_u$ being certain sheaves 
on $\blk{Y}$ with global sections $A_u$. 
Using the above diagram, 
we see
$$
\Gamma(Y,\mathcal{O}(\vartheta^*(\funcD(u))))
\; = \; 
\frac{1}{s(u)} \cdot A_u
\; \subset \: 
\K(Y).
$$
On the other hand, our present $X$ comes from the
pp-divisor $\poldiva$.
Thus, there is a canonical multiplicative 
map, forgetting the grading: 
$$ 
\bigcup_{u \in \omega \cap M} \Gamma(Y,\poldiva(u))
\; \to \; 
\K(Y),
\qquad 
f_u \mapsto f_u.
$$
This map extends to the multiplicative system of 
all homogeneous rational functions on $X$,
and hence we may may view $s(u)$ as an element of $\K(Y)$.
This gives
$$
\Gamma(Y,\mathcal{O}(\vartheta^*\funcD(u)))
=
\frac{1}{s(u)} \cdot A_u 
=
\frac{1}{s(u)} \cdot \Gamma(Y,\mathcal{O}(\poldiva(u)))
=
\Gamma(Y,\mathcal{O}(\poldiva(u)-\ddiv(s(u)))).
$$
By Lemma~\ref{compdiv}, this implies 
that $\vartheta^*(\blk{\poldiva}(u))$
equals $\poldiva(u)-\ddiv(s(u))$ for every 
$u \in\gewkegel\cap\dualgittera$.
It follows that $u \mapsto s(u)$ defines the desired
plurifunction.
\end{proof}

\setcounter{step}{0}

\begin{proof}[Proof of Theorem~\ref{mainab}]
Writing
$X:=\funcX(\poldiva)$ and
$X':=\funcX(\poldiva')$, we are given  a dominant, equivariant 
morphism $\varphi \colon X \to X'$.
The ring $A := \Gamma(X,\mathcal{O})$
is graded by the character lattice 
$\dualgittera$ 
of the torus $T = \Spec(\C[\dualgittera])$.
We will use the analogous notation 
$A'$, $\dualgittera'$ etc. for the $X'$-world.

Let us first consider the case that the pp-divisors
$\poldiva$ and $\poldiva'$ are GIT-constructed.
Let $\linabba^{*} \colon \dualgittera' \to \dualgittera$
denote the lattice homomorphism arising from the 
accompanying homomorphism $\til{\morphX} \colon T \to T'$.
By dominance of $\morphX \colon X \to X'$,
every element $u \in \dualgittera'$ gives rise
to a nonempty set
$$ 
\morphX^{-1}((X')^{ss}(u)) 
\; = \;
\bigcup_{f \in A_{\N u}} X_{\morphX^{*}(f)}
\; \subset \;
X^{ss}(\linabba^{*}(u)).
$$
By the construction of GIT-quotients, the 
set in the middle, and thus that on the left 
hand side, is a full inverse image under 
the quotient map 
$X^{ss}(F^{*}(u)) \to Y_{F^{*}(u)}$.
Hence, we obtain commutative diagrams 
(of dominant morphisms):
$$
\xymatrix{
X^{ss}(F^{*}(u)) 
\ar[d]_{\quot T}
&
{\morphX}^{-1}((X')^{ss}(u))
\ar@{}[l] |{\supset}
\ar[r]^(.6){\morphX}
\ar[d]_{\quot T}
&
(X')^{ss}(u)
\ar[d]^{\quot T'}
\\
Y_{F^{*}(u)}
&
{\morphX}^{-1}((X')^{ss}(u)) \quot T
\ar@{}[l] |{\supset}
\ar[r]_(.65){\psi_{u}}
&
Y'_{u}
}
$$

Now we consider the normalizations $Y$ and $Y'$
of the canconical components of the respective 
limits of the GIT-quotients of  $X$ and $X'$. 
The above $\psi_{u}$ fit together to a 
dominant rational map $Y \dasharrow Y'$,
defined over some open $V \subset Y$,
and we have a commutative diagram
$$
\xymatrix{
& 
V \ar[dl] \ar[dr]
&
\\
Y \ar@{-->}[rr] \ar[d]
& & 
Y' \ar[d]
\\
Y_0 \ar[rr]
& &
Y'_0}
$$
Let $\til{Y}$ denote the 
normalization of the closure 
of the graph of $V \to Y'$
in $Y \times Y'$.
Then, with the projections
$\kappa \colon \til{Y} \to Y$ and 
$\psi \colon \til{Y} \to Y'$,
we obtain a new commutative 
diagram
$$
\xymatrix{
& 
{\til{Y}} \ar[dl]_{\kappa} \ar[dr]^{\psi}
&
\\
Y \ar@{-->}[rr] \ar[d]
& & 
Y' \ar[d]
\\
Y_0 \ar[rr]
& &
Y'_0
}
$$
Since $\til{Y}$ is projective over 
the graph of $Y_0 \to Y'_0$ and hence 
over $Y_0$, 
the birational map 
$\kappa \colon \til{Y} \to Y$ 
is also projective.
Moreover, under the identification 
$\K(X)^T = \K(Y)$ the pullback
homomorphisms $\varphi^*$ and 
$\psi^*$ coincide.

Now, let 
$s \colon \dualgittera \to \K(X)^*$
and  
$s' \colon \dualgittera' \to \K(X')^{*}$, 
be the sections defining the 
minimal pp-diviors $\poldiva$ 
and $\poldiva'$ respectively,
compare Section~\ref{limit}.
Then we obtain a commutative diagram:

$$
\xymatrix{
A_{F^{*}(u)}
& &  
A'_{u}
\ar[ll]_{\varphi^{*}}
\\
{\Gamma}(\til{Y}, \mathcal{O}(\kappa^*\poldiva(F^{*}(u))))
\ar[u]^{\cdot s(F^{*}(u))}
& &
{\Gamma}(Y', \mathcal{O}(\poldiva'(u)))
\ar[u]_{\cdot s'(u)}
\ar[dl]^{\psi^{*}}
\\
&
{\Gamma}(\til{Y}, \mathcal{O}(\psi^{*}(\poldiva'(F^{*}(u))))
\ar[ul]^{\cdot \frac{\varphi^*(s'(u))}{s(F^*(u))}}
&
}
$$

The assignment
$u \mapsto \varphi^*(s'(u)) / s(F^*(u))$ 
defines a plurifunction
$\polfunca \in \K(Y,\gittera')$.
Using Lemma~\ref{compdiv}, 
one directly verifies that the triple 
$(\morphY,\linabba,\polfunca)$
describes a map of pp-divisors
$\kappa^*\poldiva \to \poldiva'$
with the properties claimed in 
the first part of the assertion.
To see the part concerning the case 
of an isomorphism
$\varphi \colon X \to X'$, note
that then no resolution
of indeterminacies $\til{Y} \to Y$ 
is needed: we can take $\kappa$ to be 
the identity, and 
$\psi \colon Y \to Y'$ is the induced 
isomorphism.

So, the assertion is proved in the case 
of GIT-constructed pp-divisors~$\poldiva$ 
and~$\poldiva'$.
In the slightly more general case
of minimal 
$\poldiva$ and $\poldiva'$,
the assertion follows immediately from 
Lemma~\ref{GITconstrpb} and the fact, that
by definition of minimality, the 
morphisms $\vartheta \colon Y \to \blk{Y}$ 
and $\vartheta' \colon Y' \to \blk{Y}'$
onto the normalized canonical components
are isomorphisms.

Now, we turn to the case that $\poldiva'$
is minimal but $\poldiva$ is not.
The part of the assertion concerning the case of 
an isomorphism $\varphi \colon X \to X'$ 
is easily settled by
using Lemma~\ref{GITconstrpb}
and the statement verified so far.

To see the first part of the assertion,
let $\poldiva_1 \in \PPDiv_\Q(Y_1,\sigma)$
be any GIT-constructed pp-divisor for $X$.  
Consider the canonical birational 
projective morphism
$\vartheta \colon Y \to Y_1$.
By Lemma~\ref{GITconstrpb}, 
the pullback $\vartheta^* \poldiva_1$ and 
$\poldiva$ differ only by the divisor of a 
plurifunction~$\polfunca$.
Moreover, by the preceding considerations, 
we have a projective, birational morphism
$\kappa_1 \colon \til{Y}_1 \to Y_1$, and a 
commutative diagram
$$ 
\xymatrix{
& 
{\funcX(\kappa_1^*\poldiva_1)} 
\ar[dl]_{\funcX(\kappa_1,\id,1)}^{\cong}  
\ar[dr]^{\funcX(\psi_1,F_1,\polfunca_1)} 
& 
\\
{\funcX(\poldiva_1)}  \ar[d]_{\cong}
& &
{\funcX(\poldiva')} \ar[d]^{=}
\\
X \ar[rr]_{\varphi} \ar[u]
& &
X' \ar[u] 
}
$$

Consider the fiber product
$Y \times_{Y_1} \til{Y}_1$.
Since all maps are birational, this 
space contains a nonempty open subset
projecting isomorphically onto open
subsets of $Y$ and $\til{Y}_1$. 
Let $\til{Y}$ be the normalization
of the closure of this subset,
and consider the canonical
projective, birational morphisms 
$\kappa \colon \til{Y} \to Y$ 
and $\vartheta_1 \colon \til{Y} \to \til{Y}_1$.
Then $\kappa^*(\poldiva)$ and 
$\vartheta_1^*\kappa_1^*\poldiva_1^*$ differ
only by the divisor of the plurifunction
$\kappa^* \polfunca$.
This allows us to define the desired map
$\kappa^*\poldiva \to \poldiva'$.

Finally, we turn to the general case.
Let $\poldiva_1' \in \PPDiv_\Q(Y_1',\sigma')$
be a GIT-constructed pp-divisor for $X'$.
Then, by what we proved so far, there
is a projective birational map 
$\kappa_1 \colon \til{Y}_1 \to Y$
and a commutative diagram.
$$ 
\xymatrix{
& 
{\funcX(\kappa_1^*\poldiva)} 
\ar[dl]_{\funcX(\kappa_1,\id,1)}^{\cong}  
\ar[dr]^{\funcX(\psi_1,F_1,\polfunca_1)} 
& 
\\
{\funcX(\poldiva)}  \ar[d]_{=}
& &
{\funcX(\poldiva'_1)} \ar[d]^{\cong}
\\
X \ar[rr]_{\varphi} \ar[u]
& &
X' \ar[u] 
}
$$

Let $\vartheta' \colon Y' \to Y_1$ be the canonical 
projective map such that $(\vartheta')^*\poldiva_1'$
and $\poldiva'$ differ only by the divisor of a 
plurifunction.
Consider the fiber product of
$\til{Y}_1$ and $Y'$ over $Y_1'$
and, similarly as before, the normalization 
$\til{Y}$ of the canonical component. 
Then we have canonical birational projective morphisms
$\kappa_2 \colon \til{Y} \to \til{Y}_1$
and $\psi \colon \til{Y} \to Y'$.
Set $\kappa := \kappa_1 \circ \kappa_2$.
Then $\kappa_1^*{\poldiva} \to \poldiva_1'$
lifts to a map 
$\kappa^* \poldiva \to (\vartheta')^*\poldiva_1'$,
which allows to define the desired map 
$\kappa^*{\poldiva} \to \poldiva'$.
\end{proof}

\section{The orbit decomposition}
\label{sec:orbits}

In this section, we use the language of 
polyhedral divisors to study the 
orbit decomposition of a normal 
affine variety with torus action.
We determine the orbit cones of Definition \ref{orbitdatadef},
and we 
describe the collection of orbits in
terms of a defining pp-divisor.
As an application, we show how 
to compute the GIT-fan of an affine 
variety with torus action directly from
its defining pp-divisor.

Let us fix the setup.
As usual, $Y$ is a semiprojective
variety, 
$N$ is a lattice with dual lattice $M$,
and 
$\sigma \subset N_\Q$ is a pointed 
cone. 
Let $\mathfrak{D} \in \PPDiv_\Q(Y,\sigma)$, 
and denote the associated
sheaf of graded algebras by
\begin{eqnarray*}
\mathcal{A}
& := & 
\bigoplus_{u \in \sigma^\vee \cap M}
\mathcal{O}(\mathfrak{D}(u)).
\end{eqnarray*}
Then we have the variety 
$\til{X} := \Spec_{Y}(\mathcal{A})$,
the torus $T := \Spec(\K[M])$, 
and the canonical map
$\pi \colon \til{X} \to Y$.
Moreover, there is a $T$-equivariant
contraction $r \colon \til{X} \to X$
onto the affine $T$-variety $X = \Spec(A)$,
where $A := \Gamma(Y,\mathcal{A})$.

Our task is to describe the $T$-orbits of 
$X$ in terms of $\poldiva$.
In Definition~\ref{fibredatadef}, we associated 
to any point $y \in Y$ a fiber polyhederon 
$\Delta_y \subset M_\Q$ with normal quasifan 
$\Lambda_y$, and a fiber monoid
complex $S_y$. 
By Corollary~\ref{orbitsfiber},
there is a bijection
\begin{eqnarray*}
\{(y,F); \; y \in Y, \; F \preceq \Delta_y\}
& \to &
\{T \text{-orbits in } \til{X} \}
\\
(y,F)
& \mapsto &
B_{\til{X}}(y,F),
\end{eqnarray*}
where $B_{\til{X}}(y,F) \subset \pi^{-1}(y)$
is the unique $T$-orbit having 
$\lambda(F) \in \Lambda_y$ as its orbit 
cone.
Besides these orbit data, our 
description of the collection of 
$T$-orbits in $X$  
involves the canonical maps
$$
\vartheta_u \colon Y \; \to \; Y_u,
\qquad
\text{where } 
Y_u \; = \; 
\Proj\left(\bigoplus_{n \in \Z_{\ge 0}} 
\Gamma(Y,\mathcal{O}(\mathfrak{D}(nu)))
\right),
$$
being induced from the semiample divisors $\mathfrak{D}(u)$. 
In a neighborhood of $y\in Y$, the
variety $Y_u$ as well as the
map $\vartheta_u$ does not depend on $u$, but only on the $\Lambda_y$-cone
containing $u$ in its relative interior.

\begin{theorem}
\label{mainresult1}
The $T$-equivariant contraction map
$r \colon \til{X} \to X$ induces a surjection
\begin{eqnarray*}
\{(y,F); \; y \in Y, \; F \preceq \Delta_y\}
& \to &
\{T \text{-orbits in } X \}
\\
(y,F)
& \mapsto &
B_X(y,F) := r(B_{\til{X}}(y,F)).
\end{eqnarray*}
One has $B_X(y_1,F_1) = B_X(y_2,F_2)$
if and only if the following two conditions
are satisfied:
$$
\lambda(F_1) 
\; = \; 
\lambda(F_2)
\; \subset \;
M_{\Q},
\qquad 
\vartheta_u(y_1) 
\; = \; 
\vartheta_u(y_2)
\text{ for some } 
u \in \relint{\lambda_{F_i}}.
$$
Moreover, for the geometry of the 
$T$-orbit $B_X(y,F) \subset X$ associated to 
a pair $(y,F)$, one obtains the following.
\begin{enumerate}
\item 
For any $x \in B_X(y,F)$, its orbit cone is
given by $\omega(x) = \lambda(F)$,
and its orbit lattice is the sublattice 
$M(x) \subset M$ generated by 
$S_y \cap \lambda(F)$.
\item 
The $T$-equivariant map of $T$-orbit 
closures
$r \colon \blk{B_{\til{X}}(y,F)} \to  \blk{B_X(y,F)}$
is the normalization.
\end{enumerate}
\end{theorem}

Note that, in contrast to the toric setting,
nonnormal orbit closures show up  quite
frequently for actions of small tori,
The simplest example is the action of $\K^*$
on $\K^2$ by means of the weights $2$ and 
$3$ --- there Neil's parabola occurs as the generic
orbit closure.

For the proof of Theorem~\ref{mainresult1}, 
we first
provide an auxiliary statement involving
the GIT-quotient 
$q_u \colon X^{ss}(u) \to Y_u$
associated to a vector $u \in M$.  

\begin{lemma}
\label{mumford}
Let $x_1,x_2 \in X$. Then we have 
$T \mal x_1 = T \mal x_2$
if and only if 
$\omega(x_1) = \omega(x_2)$ and
$q_u(x_1) = q_u(x_2)$ for some
$u \in \relint{\omega(x_i)}$.
\end{lemma}

\begin{proof}
Only the ``if'' part is nontrivial.
So, suppose $\omega(x_1) = \omega(x_2)$ and
$q_u(x_1) = q_u(x_2)$ for some
$u \in \relint{\omega(x_i)}$.
According to Proposition~\ref{orbitcones},
any $f \in A_{nu}$, where $n > 0$,
vanishes along
$\blk{T \mal x_i} \setminus T \mal x_i$.
Consequently, the $T$-orbits through $x_1$ 
and $x_2$ are closed in $X_u$.
Since good quotients separate closed orbits,
$q_u(x_1) = q_u(x_2)$ implies
$T \mal x_1 = T \mal x_2$.
\end{proof}

\begin{proof}[Proof of Theorem~\ref{mainresult1}]
We first prove statements~(i) and~(ii)
on the geometry of the orbits.
Fix a pair $(y,F)$,
and choose a point 
$\til{x} \in \pi^{-1}(y)$ with 
orbit cone $\omega(\til{x}) = \lambda(F)$. 
Then the associated orbit lattice
$M(\til{x}) \subset M$ is generated
by $S_y \cap \lambda(F)$.
We will show~(i) by checking that 
$x := r(\til{x})$
has orbit data
$\omega(x) = \omega(\til{x})$ and
$M(x) = M(\til{x})$.
Since the $T$-orbit closure of $\til{x}$
is normal, this also proves~(ii).
 
In order to see $\omega(x) \subset \omega(\til{x})$,
let $u \in \omega(x)$. 
Then there is an $f \in \Gamma(X,\mathcal{O})_{nu}$
with $n>0$ such that $f(x) \ne 0$
holds. Thus, $r^*f(\til{x}) \ne 0$, which
implies $u \in \omega(\til{x})$.
For the reverse inclusion, note
that we find for every $u \in \omega(\til{x})$
a $g \in \Gamma(Y,\mathcal{A}_{nu})$,
where $n > 0$, such that 
$\pi(\til{x}) \not\in Z(g)$.
Then we have $g(\til{x}) \ne 0$.
Moreover, $g = r^*f$ with 
$f \in \Gamma(X,\mathcal{O})_{nu}$
and $f(x) \ne 0$.
This implies $u \in \omega(x)$.

Similar to the orbit cones, we see
$M(x) \subset M(\til{x})$. 
To verify the reverse inclusion, 
let $u \in  S(\til{x})$.
Consider the contraction map
$\vartheta_u \colon Y \to Y_u$.
Then, $\mathfrak{D}(u) = \vartheta_u^*(E_u)$
with an ample divisor $E_u$ on $Y_u$.
In particular, we have 
$\mathfrak{D}(u) = \kdiv(h^{-1})$ on some 
neighbourhood $V = \vartheta_u^{-1}(V_u)$
of $y = \pi(\til{x})$ with $V_u \subset Y_u$ open
and $h \in \K(V)$. 

Recall that there is a good quotient
$X_u \to Y_u$ for the set 
$X_u \subset X$ 
of semistable points associated to
$u \in \omega \cap M$. Moreover, 
we may restrict $\pi \colon \til{X} \to Y$
to obtain a morphism
$r^{-1}(X_u) \to Y$.
Denoting by $W \subset r^{-1}(X_u)$ 
and $W_u \subset X_u$ the inverse images 
of $V \subset Y$ and $V_u \subset Y_u$
respectively, we arrive at a commutative 
cube:
$$
\xymatrix@!0{
& 
{r^{-1}(X_u)} \ar[rr]  \ar'[d][dd]
& & 
X_u \ar[dd] 
\\
W \ar[rr] \ar[ur]  \ar[dd]
& & 
W_u  \ar[ur] \ar[dd]
\\
&
Y \ar'[r][rr] 
& & 
Y_u 
\\
V \ar[rr]  \ar[ur] 
& & 
V_u  \ar[ur]
}
$$

Now, consider $h \in \K(V)$ as a 
regular function on 
$W \subset \pi^{-1}(V)$. 
Since $u \in S(\til{x}) \subset \omega(\til{x})$ 
and $\omega(\til{x}) = \omega(x)$,
we have $x \in X_u$.Hence $\til{x} \in r^{-1}(X_u)$,
which gives us $\til{x} \in W$. 
Moreover, since $y \not\in Z(h)$ holds, 
$h$ is not trivial along $\pi^{-1}(y)$,
and thus $u \in S(\til{x})$ yields 
$h(\til{x}) \ne 0$. 
Since $W \to W_u$ is proper and birational,
the function $h \in \Gamma(W,\mathcal{O})_u$
is in fact a regular function on $W_u$.
By construction, we have $x \in W_u$ and
$h(x) \ne 0$.
This implies $u \in M(x)$.

We come to the characterization
of the equality 
$B_X(y_1,F_1) = B_X(y_2,F_2)$.
Choose $x_i \in B_X(y_i,F_i)$.
As we have just seen,
$\omega(x_i) = \lambda(F_i)$
holds.
Moreover, as remarked just 
before,
we have a commutative diagram 
for every $u \in \omega \cap M$ 
$$
\xymatrix{
r^{-1}(X_u) \ar[r]^{r} \ar[d]_{\pi} 
&
X_u \ar[d]^{q_u}
\\
Y \ar[r]_{\vartheta_u}
&
Y_u.
}
$$
Thus, the conditions
$\lambda(F_1) = \lambda(F_2)$ 
and $\vartheta_u(y_1) = \vartheta_u(y_2)$
are equivalent to the conditions
$\omega(x_1) = \omega(x_2)$
and $q_u(x_1) = q_u(x_2)$.
According to Lemma~\ref{mumford},
the latter conditions characterize 
$T \mal x_1 = T \mal x_2$.
\end{proof}

Putting together Theorem~\ref{mainresult1}
with 
Corollary~\ref{orbitsfiber}
and Propositions~\ref{pifibers}, ~\ref{reduced}
gives the following characterization
for a pp-divisor to be an integral Cartier 
divisor.

\begin{corollary}
The following statements are equivalent.
\begin{enumerate}
\item 
The polyhedral divisor $\poldiva$  
belongs to $\CDiv(Y,\sigma)$.
\item 
The map $\pi \colon \til{X} \to Y$
has no multiple fibers.
\item 
The torus $T$ acts with connected isotropy 
groups on $X$.
\end{enumerate}
\end{corollary}

In a further application, we indicate
how to read off the GIT-quasifan of $X$ 
in the sense of Theorem~\ref{gitfan}
from its defining pp-divisor 
$\poldiva$.

\begin{corollary}
Let 
$
\mathfrak{D}  
=  
\Delta_1 \otimes D_1+ \ldots +   \Delta_r \otimes D_r
$
with prime divisors $D_i$.
Then the quasifan of GIT-cones associated 
to the $T$-linearizations of the
trivial bundle on $X$ is the normal 
quasifan of the Minkowski sum 
$\Delta_1 + \ldots + \Delta_r$.
\end{corollary}

\begin{proof}
According to Theorem~\ref{gitfan}, 
the GIT-quasifan $\Lambda$ of $X$ 
is the coarsest quasifan in $M_\Q$ 
refining 
all orbit cones $\omega(x)$, 
where $x \in X$.
By Theorem~\ref{mainresult1}, the orbit 
cones $\omega(x)$, where $x \in X$,
are precisely the cones of the 
normal quasifans $\Lambda_y$ 
of the fiber polyhedra
$\Delta_y$, where $y \in Y$.
Thus, $\Lambda$ is the coarsest
common refinement of all $\Lambda_y$ 
and hence equals the normal quasifan 
of the Minkowski sum of all the 
$\Delta_y$ where $y \in Y$.
But the latter equals the normal 
quasifan of
$\Delta_1 + \ldots + \Delta_r$.
\end{proof}

\section{Calculating Examples}
\label{sec:examples}

In this section, we indicate a recipe
how to determine a minimal pp-divisor
for a given normal affine $T$-variety
$X$. The strategy is first to treat
the case of a toric variety $X$ and 
then to settle the general case via
equivariant embedding.
The proof of the method is straightforward
and will therefore be ommitted.

Consider an affine toric 
variety $X$ and the action of a subtorus 
$T \subset T_X$ of the big torus
$T_X \subset X$.
Let $N_X$ be the lattice of 
one parameter subgroups of $T_X$,
and let $\delta \subset (N_X)_\Q$ 
be the cone describing $X$. 
The inclusion $T \subset T_X$ 
corresponds to an inclusion
$N=N_T \subset N_X$ of lattices,
and we obtain a (non-canonical) split 
exact sequence
$$
\xymatrix@1{
0 \ar[r] &
N_T \ar[r]^{F} & 
N_X \ar[r]^{P} \ar@/^1pc/[l]_-{s} &
N_Y \ar[r] &
0,}
$$
where $N_Y := N_X/N_T$ and $s \colon N_X \to N_T$
satisfies $s \circ F = \id$.
Let $\Sigma_Y$ be the coarsest fan in 
$(N_Y)_\Q$ refining all cones $P(\delta_0)$
where $\delta_0 \preceq \delta$.
Then the toric variety $Y$ corresponding 
to $\Sigma_Y$ is the normalization of the 
closure of the image of $T_X$ in the limit 
over all GIT-quotients of $X$,
i.e.\ $Y$ is as in
Section~\ref{limit}.
Note that, up to normalization, $Y$ equals the Chow quotient of $X$ by $T$
as constructed in \cite{Chow}.

Let us indicate how to obtain the minimal 
pp-divisor for the $T$-variety $X$. 
Given a one-dimensional 
cone $\varrho \in \Sigma_Y$,
let $v_\varrho \in \varrho$
denote the first lattice vector, 
and define a polyhedron
$$ 
\Delta_{\varrho}
\; := \;
s(\delta \cap P^{-1}(v_\varrho))
\; \subset N_\Q=(N_T)_\Q.
$$
These polyhedra have
$\sigma := \delta \cap (N_T)_\Q$
as their tail cone.
Denoting by $\Rho_Y \subset \Sigma_Y$
the set of one-dimensional cones,
and by 
$D_\varrho \subset Y$ 
the prime divisor corresponding 
$\varrho \in \Rho_Y$, 
we obtain a minimal pp-divisor
for the $T$-variety $X$ as    
$$ 
\mathfrak{D}
\; = \; 
\sum_{\varrho \in \Rho_Y} \Delta_{\varrho} \otimes D_{\varrho}
\; \in \; 
\PPDiv(Y,\sigma).
$$   

\begin{example}
Consider the affine toric variety $X = \K^4$ and the 
action of the two-dimensional torus $T = (\K^*)^2$
on $X$ given, with respect to standard coordinates, 
by
\begin{eqnarray*}
t \mal z 
& = & 
(t_1^4z_1, t_1^3z_2, t_2z_3, t_1^{12}t_2^{-1}z_4).
\end{eqnarray*}
The corresponding lattices are $N_T = \Z^2$,
and $N_X = \Z^4$, and the quotient lattice
is $N_Y = \Z^2$. The maps $F \colon N_T \to N_X$ 
and $P \colon N_X \to N_Y$ and 
a choice for $s \colon N_X \to N_T$ 
are given by the 
matrices
$$ 
F \; = \; 
\left[
\begin{array}{rr}
4  & 0 \\
3  & 0 \\
0  & 1 \\
12 & -1
\end{array}
\right],
\qquad
s \; = \; 
\left[
\begin{array}{rrrr}
1  & -1 & 0 & 0 \\
0  & 0  & 1 & 0
\end{array}
\right],
\qquad
P
\; = \; 
\left[
\begin{array}{rrrr}
3 & 0 & -1 & -1 \\
0 & 4 & -1 & -1
\end{array}
\right].
$$
>From this we see that the fan 
$\Sigma_Y$ is the standard fan of 
the projective plane, meaning that it
has as its maximal cones
$$ 
\cone((1,0),(0,1)),
\qquad
\cone((0,1),(-1,-1)),
\qquad
\cone((-1,-1),(1,0)).
$$
The cone $\sigma = s(\Q^4_{\ge 0} \cap F(\Q^2))$
is generated by 
the vectors $(1,0)$ and $(1,12)$.
The polyhedral coefficients of the 
pp-divisor $\sum \Delta_\varrho \otimes D_\varrho$ 
on $Y = \PP^2$ are
$$
\begin{array}{lclcl}
\Delta_{\Q_{\ge 0}(1,0)} 
& = & 
\left(1/3,0\right) \; + \; \sigma
& = : &
\Delta_0,
\\
\Delta_{\Q_{\ge 0}(0,1)} 
& = & 
\left(-1/4,0\right) \; + \; \sigma
& = : & 
\Delta_1,
\\
\Delta_{\Q_{\ge 0}(-1,-1)} 
& = & 
\left(\{0\}\times[0,1]\right) \; + \; \sigma
& =: & 
\Delta_\infty.
\end{array}
$$
\end{example}

Let us indicate how to handle the general 
case, i.e.\ a possibly non-toric normal 
affine variety $X$ with an effective action 
of a torus $T$.
First, choose a $T$-equivariant embedding 
into some $\K^n$, where $T$ acts as 
a subtorus of $(\K^*)^n$, and 
$X$ hits the big orbit of
$(\K^*)^n$.
Then apply the previous method 
to obtain a minimal pp-divisor
\begin{eqnarray*}
\poldiva_{\rm toric}
& = &
\sum \Delta_\varrho \otimes D_\varrho
\end{eqnarray*}
for the $T$-variety $\K^n$ living 
on some toric variety $Y_{\rm toric}$.
Then the desired variety $Y$ lying 
over the GIT-quotients of $X$ is 
the normalization of the closure 
of the image of $X \cap (\K^*)^n$
in $Y_{\rm toric}$.
Moreover, a minimal pp-divisor for
the $T$-variety $X$ is obtained by 
pulling back $\mathfrak{D}_{\rm toric}$ 
to $Y$.

\begin{example}
Consider once more the affine threefold
$X = V(z_1^3+z_2^4+z_3z_4)$ in $\K^4$
discussed in the introduction of the 
paper.
It is invariant under the action
of $T = (\K^*)^2$ on $\K^4$ given by 
\begin{eqnarray*}
t \mal z 
& = & 
(t_1^4z_1, t_1^3z_2, t_2z_3, t_1^{12}t_2^{-1}z_4).
\end{eqnarray*}
In the preceding example, we computed 
a minimal pp-divisor $\mathfrak{D}_{\rm toric}$
living on $Y_{\rm toric} = \PP^2$ for
the $T$-action on $\K^4$.
On the big tori $(\K^*)^4 \subset \K^4$ 
and $(\K^*)^2 \subset \PP^2$,
the projection is given by
$$
(\K^*)^4 \; \to \; (\K^*)^2,
\qquad
(t_1,t_2,t_3,t_4) \; \mapsto \;
\left(\frac{t_1^3}{t_3t_4},\frac{t_2^4}{t_3t_4}\right).
$$
>From this we see that the closure of the 
image of $X \cap (\K^*)^4$ 
in $\PP^2$ is given with respect to its
homogeneous coordinates $(w_0:w_1:w_2)$ by 
\begin{eqnarray*}
Y 
& = & 
V(w_0+w_1+w_2).
\end{eqnarray*}
Since $Y = \PP^1$ is already normal, 
we can restrict $\mathfrak{D}_{\rm toric}$ 
to $Y$ and thus obtain a minimal pp-divisor
for the $T$-variety $X$ by 
\begin{eqnarray*}
\poldiva
& = &
\Delta_0 \otimes \{0\} 
\; + \; 
\Delta_1 \otimes \{1\}
\; + \;
\Delta_\infty \otimes \{\infty\}. 
\end{eqnarray*}
\end{example}

\end{document}